\documentclass[10pt]{article}

\usepackage{float}
\usepackage{amssymb, color, epsfig, graphics,mathtools}
\usepackage[dvipsnames]{xcolor}
\usepackage{graphicx}
\usepackage{subcaption}
\usepackage{amsmath}
\usepackage{latexsym} 
\usepackage[all]{xy}
\usepackage{hyperref}

\newtheorem{proposition}{Proposition}

\numberwithin{equation}{section}

\usepackage{geometry,multicol}
\newdimen\stockheight
\newdimen\stockwidth
\stockwidth\paperwidth
\newcommand{\be}{\begin{equation}}
\newcommand{\ee}{\end{equation}}
\newcommand{\bea}{\begin{eqnarray}}
\newcommand{\eea}{\end{eqnarray}}
\newcommand{\barr}{\begin{array}}
\newcommand{\earr}{\end{array}}
\newtheorem{rem}{Remark}[section]

\newcommand{\bpar}{\begin{equation} \left\{ \begin{array}{lll}}
\newcommand{\epar}{ \end{array}\right. \end{equation} }
\newcommand{\eparn}{ \end{array} \right.}

\newcommand\eqs{equations }

\newcommand\leftmat{\left(\begin{array}{cc}}
\newcommand\rightmat{\end{array}\right)}
\newcommand\leftvec{\left(\begin{array}{c}}
\newcommand\rightvec{\end{array}\right)}
\newcommand\re{{\rm e}}

\newcommand\R{{\mathbb R}}
\newcommand\C{{\mathbb C}}
\newcommand\Z{{\mathbb Z}}

\newcommand{\Li}{\operatorname{Li}} 
\newcommand{\Cl}{\operatorname{Cl}} 
\newcommand{\Sl}{\operatorname{Sl}} 
\newcommand{\sign}{\operatorname{sgn}} 

\newcommand{\dide}{\delta_{[-\pi,\pi]}}

\input otex.sty

\def\i{\:\ro i\:}

\def\Ft{\CF}
\def\Ftr#1{\Ft[\:#1\:]}
\def\Ftrx#1{\Ftr{#1(x)}}

\def\i{\:\ro i\:}
\def\L{\Bbb L}

\def\Id{{\cal I}_\delta}
\def\Sd{{\cal S}_\delta}

\begin{document}

\title{New Revival Phenomena for \break Linear Integro--Differential Equations}
\author{
Lyonell Boulton\textsuperscript{\textasteriskcentered}, Peter J. Olver\textsuperscript{\textdagger}, Beatrice Pelloni\textsuperscript{\textasteriskcentered}, David A. Smith\textsuperscript{\textdaggerdbl} \\
{\small \textsuperscript{\textasteriskcentered} Department of Mathematics, Heriot-Watt University, UK} \\ {\small \textsuperscript{\textdagger} School of Mathematics, University of Minnesota, MN, USA} \\ {\small \textsuperscript{\textdaggerdbl} Corresponding author \href{mailto:dave.smith@yale-nus.edu.sg}{\texttt{dave.smith@yale-nus.edu.sg}} } \\ {\small Division of Science, Yale-NUS College, Singapore}
}
\date{1st October 2020}
\maketitle

\begin{abstract}
    We present and analyse a novel manifestation of the revival phenomenon  for linear spatially periodic evolution equations, in the concrete case of three nonlocal equations that arise in water wave theory and are defined by convolution kernels.
    Revival in these cases is manifested in the form of dispersively quantised cusped solutions at rational times.
    We give an analytic description of this phenomenon, and present illustrative numerical simulations.
\end{abstract}

\section{Introduction}

\emph{Revival} for linear partial differential equations comes from the experimentally observed phenomenon of \emph{quantum revival}~\cite{BMS, VVS, YeaStr}, in which an electron that is initially concentrated near a single location of its orbital shell is re-concentrated near a finite number of orbital locations at certain times.  A precursor of  this phenomenon was observed as far back as 1834 in optical experiments performed by Talbot \cite{Talbot}.  This motivated the pioneering work of Berry and collaborators \cite{Berry, BerryKlein, BMS}, on what they called the Talbot effect in the context of the linear free space Schr\"odinger equation. We refer the reader to the introductory text \cite{Thaller} for an account on this.  The concept was extended later to a class of linear dispersive equations that included the linearised Korteweg--deVries equation, first by Oskolkov, \rf{OskolkovV} and subsequently rediscovered and extended by one of the authors \cite{Odq}, who called the effect \emph{dispersive quantisation}. See the recent survey~\cite{Smi2020a}.

As a result of these investigations, one learns that a dispersive linear partial differential equation in one space variable, possessing a polynomial dispersion relation and subject to periodic boundary conditions, exhibits a  \emph{revival phenomenon}.   The corresponding fundamental solution, i.e the solution evolving from a delta function (extended periodically) at time zero, at certain values of time localises into a suitable finite linear combination of copies of this delta function. Specifically, this occurs for a dense set of times that includes all rational multiples of a fixed quantity depending on the equation and the period. This has the striking consequence that the solution to any initial value problem, at such \emph{rational times}, depends only upon finitely many values of the initial datum. The latter constitutes the simplest mathematical manifestation of the revival phenomenon and corresponds precisely to the experimental observations mentioned above. 

If the initial condition at time zero is piecewise smooth, at such rational times the solution, being a linear combination of piecewise smooth translates of this initial condition, is piecewise smooth. Similarly, if the initial condition is piecewise constant, the resulting solution is also piecewise constant at rational times. But, remarkably, if the initial condition is discontinuous, the solution profile at the other ``irrational'' times turns out to be a continuous nowhere differentiable function with graph of a specific fractal dimension \cite{ErSh,ET}. The latter effect is known as \emph{dispersive fractalisation}. Strong evidence, backed by numerical computations, suggests that such fractalisation occurs for all cases when the revival property holds.

It is worth noting here that the existence of revival and fractalisation also depends on specific boundary conditions, although in a manner not fully understood at present. Promising preliminary results in this direction were reported in \rf{OSS}.
     
Following a different line of enquiry, in \cite{COdisp}, the phenomena of fractalisation and dispersive quantisation were numerically shown to extend to nonlinear dispersive wave equations on periodic domains, including the integrable nonlinear Schr\"odinger  and Korteweg--deVries equations as well as their non-integrable generalisations. For these models, the solution at rational times assumes a quantised form that is noticeably different from its nearby fully fractal profiles. Roughly, it is much smoother, albeit discontinuous, or at least less markedly fractal. While some rigorous results on the fractilisation at a dense set of time have been established, \rf{ChErTz}, for the most part these are only  numerical observations, valid for small times, but they indicate a tantalising robustness of the revival phenomenon. Moreover, following an implementation of a state-of-the-art numerical scheme for the solution of nonlinear time evolution problems \rf{KS, KS1}, this has recently been confirmed for the nonlinear Schr\"odinger equation also in the large time regime. 
     
Several other examples of linear integro-differential equations with non-polynomial dispersion relations, mostly arising from fluid mechanics (water waves) and again subject to periodic boundary conditions, were also investigated in \cite{COdisp}. A unit step, also known as the Riemann problem \rf{Whitham}, was considered as initial condition. A wide range of phenomena with qualitative features depending upon the large wave number asymptotics of the dispersion relations were observed.  Among such equations, those that have an asymptotically polynomial dispersion relation, which include the linearised Benjamin-Ono and integrable Boussinesq equations, were shown to possess dispersively quantised cusped solutions at rational times. Interestingly, these solutions appear to be piecewise smooth but non-constant between cusps.

The concrete goal of the present paper is to investigate the phenomenon of dispersively quantised cusped solutions in greater detail, by examining the periodic Riemann problem for three linear integro-differential equations that arise in models from fluid mechanics. We consider the  Benjamin--Ono (BO) equation, \rf{ABFS,  Benjamin,Ono,DP, Saut}, the Intermediate Long Wave (ILW) equation, \rf{ABFS,  Joseph,KKD,DP, Saut}, and a unidirectional version of the integrable Boussinesq equation, \rf{MCK}, known as the Smith\footnote{The Smith for which this equation is named is no relation to the author of the present work.} equation, \rf{ABFS, Smith}. These three equations are intimately related. As we describe below, the BO equation is obtained as the  limit of the ILW equation and of the Smith equation when the parameter modelling water depth tends to infinity.

We provide an explicit characterisation of the previously observed cusped profiles, leading to a more general form of revival for linear integro-differential equations with non-polynomial dispersion relations. Clearly, this new revival phenomenon is worth further study, both mathematically and in terms of how it may impact applications.

Our main result is an analytic description of this new form of revival. For the BO equation we show that the solution of the periodic Riemann problem is a {\em linear combination of translates of the initial step function and of its periodic Hilbert transform}.  For the ILW and Smith equations, a similar result holds but only in an approximate sense.  For these two equations, we provide a quantification of the error.


In the course of our analysis, we will introduce a certain class of special functions, that we call {\em trigonometric polylogarithms}. A particularly important subclass is given by the {\em trigonometric hypergeometric functions}, so named in that they are obtained by evaluating classical hypergeometric functions on the unit circle. The latter functions can, in fact, be expressed in terms of elementary functions. Various properties of these special functions play a central role in expressing the new type of revival phenomena observed in the above linear integro-differential equations. 

\section{Three Linearised Model Wave Equations} \label{sec:ModelWaveEquations}

The unidirectional evolution equations we consider can be expressed in terms of their dispersion relation $\omega (k)$ connecting the wave number $k$ (or spatial frequency) to the temporal frequency $ \omega $, \rf{Whitham}. Indeed, the equations are of the general form 
\Eq{leq}
$$
u_t=L[u],$$
in which $L$ is a linear integro-differential operator characterised by its Fourier transform 
\Eq{leqF}
$$\widehat{Lu}(k)=\omega (k)\,\widehat u(k) .
$$
We use the Fourier transform normalisation
\Eq{Ftr}
$$\fhat(k) =\Ftrx f = \fra{\sqrt {2\pii }}\> \intiix f(x)\,{e^{-\i k\:x}} ,$$
so that the Fourier transform of the convolution of two functions,
$$f * g\,(x) = \intiiy f(x-y)\,g(y) $$
is given by
$$ \widehat{f*g}\,(k)=\sqrt {2\pii } \  \fhat(k) \, \ghat(k).$$

\subsection{The Linearised Benjamin--Ono Equation}
\label{lBO}

 The  Benjamin--Ono equation  arises as a nonlinear integrable model for internal waves in fluids, and has been studied extensively,  \rf{ABFS, Benjamin,Ono,Saut}.  Its linearisation has the form
\Eq{lBO}
$$u_t = \CH[u_{xx}],$$
where $\CH$ denotes the \is{Hilbert transform}
\Eq{Htr} 
$$\CH[f](x) = H * f\,(x) = \frac 1{\pi} \ \intii y{\hskip-24pt \hbox{---} \hskip14pt {\frac{f(y)}{x-y}}} ,$$
for the Cauchy kernel
$H(x) = 1/\pa{\pi\, x}$.  The bar across the integral sign denotes principal value.
We are interested in the spatially periodic problem on the interval $\ipp$.  In this case, the convolution formula \eq{Htr} is replaced by the \is{circular} or \is{periodic Hilbert transform}, \rf{King}, \cite[Ch.~9]{BuNe}:
\be{}\CH[f](x) =\frac 1{\pi} \Sumii k\, \intpp y{\hskip-21pt \hbox{---} \hskip10pt {\frac{f(y)}{x-y+ 2\pii\:k}}} 
= \fra{2\pii}\, \intpp y{\hskip-21pt \hbox{---} \hskip10pt \cot\bbk{\f2\pa{x-y}} f(y)} ,
\label{perhtr}\ee
which is a periodic convolution with the rescaled cotangent kernel.
The dispersion relation for the linearised BO equation is the odd function 
\Eq{BOdr}
$$\omega_{BO}(k) = k^2 \sign k.$$

\subsection{The Linearised Intermediate Long Wave Equation}
\label{lILW}

The nonlinear Intermediate Long Wave  equation is yet another integrable equation that arises through the modelling of internal waves in a fluid of finite depth, \rf{ABFS,Joseph,KKD, Saut}.  Its infinite depth limit is the BO equation, as we explain below.

The linearised ILW equation is given by
\Eq{lILW}
$$u_t={\cal L}[u],$$
with the integro-differential operator ${\cal L}$ given
by
\be
{\cal L}[u]={\cal I_\delta}[u_{xx}]-\frac 1 \delta \,u_x,
\label{calL}
\ee
where
\Eq{Idconv}
$$\Id[f](x) = -\fra{2\,\delta } \ \intii y{\hskip-24pt \hbox{---} \hskip14pt \coth \Bk{\frac \pi{2\,\delta }\, (x-y)}\,f(y)} $$
is convolution with a hyperbolic cotangent kernel. If $f(x)$ is $2\pii$-periodic, this convolution formally becomes
\Eq{Idperiodic}
$$\Id[f](x) = -\fra{2\,\delta }\ \intpp y{\hskip-21pt \hbox{---} \hskip8pt \Bk{\Sumii n \coth \Pa{\frac \pi{2\,\delta }\, (x-y) + \frac{\pi^2 n}{\delta }}}\,f(y)} \ .$$
We now seek a closed form expression for a convolution kernel $C_\delta(x)$, such that 
\Eq{ConvCdelta}
$$\Id[f](x) = C_\delta  * f\,(x)=\intpp y{\hskip-21pt \hbox{---} \hskip8pt  C_{\delta}(x-y) f(y)} \ .$$ 
As for the case of \eqref{perhtr}, $C_{\delta}(x)$ will have a singularity of order $1/x$ (a simple pole, when considered as an analytic function for complex $x$) at $x=y$, so the integral should be taken as a principal value. 

As a starting point, consider the formal expression
\Eq{ILWker}
$$C_\delta (x) = -\fra{2\,\delta } \Sumii n \coth \Pa{\frac \pi{2\,\delta }\,x + \frac{\pi^2 n}{\delta }}.$$ Let $\zeta (z)$ denote the Weierstrass zeta function associated with the lattice $\L \subset \C$ generated by $2\:\omega _1, 2\:\omega _3 \in \C$. Then,  \rf{DLMF; 23.8.4}
\Eq{zetacot}
$$\xeq{\zeta (z) = \frac{\eta _1}{\omega _1}\,z +  \frac{\pi}{2\,\omega _1} 
\Sumii n \cot \Pa{\frac{\pi}{2\,\omega _1}\,z + \frac{\pi\,\omega _3\, n}{\omega _1}},\\z \not\in \L,\wherex \eta _1 = \zeta (\omega _1).}$$
 Set $\myeq{\omega _1 = -\i \delta$ and $\omega _3 = \pi}$. Then $\Im (\omega _3/\omega _1) = \pi /\delta  >0$ and $\L$ becomes a real rectangular lattice. Since $\coth z = \i \cot(\i z)$, comparing \eqs{ILWker}{zetacot}, yields
\Eq{Czeta}
$$C_\delta (x) = \fra{\pii} \bbk{\alpha \,x - \zeta (x) },\where \alpha = \i\> \frac{\zeta (-\i \delta )}{\delta }$$
for all real $x\not= 2n\omega_3=2n\pi$.  

From the periodicity properties of $\zeta(z)$ \rf{DLMF;23.2(iii)}, it follows that $C_{\delta}(x)$ is quasi-periodic with 
$$C_\delta(x+2\pi)=C_{\delta}(x)-\frac{1}{\delta}$$
and in fact periodic with period $2\i/{\delta}$. To rigorously prove that the expresion for $\Id[f](x)$ in \eq{Idconv} for $2\pi$-periodic $f(x)$ and the circular convolution \eq{ConvCdelta} coincide, it might be possible to follow a similar programme as described in \cite[\S 9.1.1]{BuNe} for the Hilbert transform. This will be pursued elsewhere.\footnote{A more complicated closed form expression for the kernel $C_\delta$ involving Jacobi elliptic functions can be found in \rf{AFSS}. The expression given here in \eq{Czeta} appears to be novel. }



The corresponding non-polynomial dispersion relation for the (linearised) ILW equation \eq{lILW}, which depends on the depth parameter $\delta > 0$, is given by the odd function
\Eq{ILWdr}
$$	  \omega _\delta(k)=k^2\coth(\delta k)-\frac k \delta.$$
Note that in the infinite depth  limit, $\delta \to \infty $, the ILW dispersion relation \eq{ILWdr} reduces to the dispersion relation \eq{BOdr} for the BO equation as $\omega _\delta(k) \longrightarrow \omega _{BO}(k)$ pointwise. 

\subsection{The Linearised Smith Equation}
\label{lsmith}

The Smith equation was proposed  in \cite{Smith} in the context of water wave theory,  as a model for continental shelf waves. An analysis of the existence and regularity of the solutions of the Cauchy and of the periodic problem for this equation is presented in \cite{ABFS}.  
The Smith equation can be viewed as the unidirectional analogue of the more famous Boussinesq equation, \rf{MCK}, whose revival properties were first investigated in \cite{COdisp}.
Indeed, the two equations have the same dispersion relation.

The Smith equation as analysed in \rf{ABFS} is given in dimensionless coordinates, and proved to be a bounded perturbation of the BO equation. To quantify this qualitative result, and to shed light on the similarity in the revival property the BO and Smith equations display, it is necessary to understand the modelling assumptions made in the derivation given in \rf{Smith}. This derivation depends on a parameter $\delta$ directly proportional to the water depth, and that we have explicitly included in the dispersion relation. The limit $\delta\to\infty$  therefore corresponds to infinite depth, namely the BO regime.

Our version of the Smith equation takes the form
\Eq{Smith}
$$u_{t}=\Sd[u_x], $$
where the integro-differential operator  is  defined through the Fourier expansion of its kernel
\Eq{Smithkernel}
$$
\shat_\delta(k) = \i\sqrt{k^2 + \fra \delta } \roq{whereby} \widehat{\Sd[u]}(k)=\shat_\delta(k)\,\uhat (k).
$$
In other words, 
\Eq{Sd}
$$\Sd[f](x)= \frac \i{2\pii} \intii k{\intii y{\sqrt{k^2+\frac 1 \delta} \;e^{\i k(x-y)}f(y) } } .$$
In the Appendix, we provide details of a computation that yields the inverse Fourier transform of \eq{Smithkernel} in closed form as
\Eq{smithif}
$$
s_\delta(x) = -\i\sqrt{\frac 2{\pi\,\delta }}\;\frac{K_1\bpa{\abs x/\sqrt \delta \,}} {\abs x},$$
where $K_1(x)$ denotes the modified Bessel function of the second kind, \rf{DLMF; 10.25}.
Thus, the Smith operator \eq{Sd} on $\mathbb R$ is given by convolution
\Eq{Smithc}
$$\Sd[f] = \frac{s_\delta *f\,(x)}{\sqrt{2\pii}}  = -\frac \i{\pi\,\sqrt \delta } \intiiy {\frac {K_1\bpa{\abs {x-y}/\sqrt \delta\,}} {\abs {x-y}}\,f(y)} .$$

If $f$ is $2\pii$-periodic, then we can write
$$\Sd[f] = S_\delta*f\,(x) = \intppy {S_\delta(x-y)\,f(y)} ,$$
where, formally,
\Eq{perSmithker}
$$S_\delta(x) = -\frac \i{\pi\,\sqrt \delta }\;\Sumii n \frac{K_1\bpa{\abs{x + 2\:n\pii}/\sqrt \delta\,}} {\abs{x + 2\:n\pii}}.$$
The Bessel function $K_1(z)$ is entire in $z$ except at $z=0$ where it has a pole of order $1$, and hence the only singularities of $S_\delta (x)$ are at the even multiples of $\pi$.  Moreover, according to \rf{DLMF; 10.40.2}, 
$$K_1(x) \simx \sqrt{\frac{\pi}{2\:x}}\ e^{-x}\roq{as} x \longrightarrow \infty ,$$ 
and hence, away from the pole singularities, the infinite sum \eq{perSmithker} converges rapidly.  However, unlike the periodic BO and ILW kernels, it seems unlikely that there exists a closed form formula for  the summation \eq{perSmithker} in terms of known special functions.

The dispersion relation for the Smith equation \eq{Smith} is
\Eq{Smithdr}
$$
    \omega_{S }(k) = k\sqrt{\frac 1\delta+k^2}.
$$
As for the case of the ILW equation, in the limit as the depth parameter becomes infinite, the Smith dispersion relation \eq{Smithdr} reduces to the dispersion relation \eq{BOdr} for the BO equation: $\omega _{S }(k) \longrightarrow \omega _{BO}(k)$ pointwise as $\delta \to \infty$.

\medskip

\Rmk Intriguingly the hierarchy (ignoring constants)
\Eq{hierarchy}
$$\fra x \longmapstoq \cot x \orx \coth x \longmapstoq \zeta (x)$$
i.e., rational to (hyperbolic) trigonometric to elliptic, corresponding to
$$\roh{BO} \longmapstox \roh{periodic BO or ILW} \longmapstox \roh{periodic ILW}$$
reminds one of the hierarchies in the discrete integrable Calogero--Moser--Sutherland many body models, both classical and quantum, \rf{Calogero, MoserH, Sutherland}.  Indeed, if one differentiates \eq{hierarchy} one obtains
$$-\,\fra{x^2} \longmapstoq -\,\csc^2 x\orx -\,\csch^2 x\longmapstoq -\,\CP(x),$$
where $\CP$ is the Weierstrass elliptic function, which are precisely the Calogero--Moser--Sutherland potentials.  Connections between these continuous and discrete integrable systems are noted in \rf{ZaZo}.


\section{Trigonometric polylogarithms}

Our analysis will make use of certain special functions that are parametrised by a positive integer
$r \in \N^+$.
In this section, we introduce these functions and present their main properties. The case $r=1$ is particularly important, hence we give a detailed account of it in this section; corresponding results for all other cases can be found in the Appendix. 


Let $j,k,r\in\N^+$, $j\leq k$.
We define the \emph{trigonometric polylogarithms}, $S^k_{j,r}(x)$ and $C^k_{j,r}(x)$, as
\Eq{SC}
$$    S^k_{j,r}(x) = \Im\left[ E^k_{j,r}(x) \right] \qquad \text{and} \qquad
    C^k_{j,r}(x) = \Re\left[ E^k_{j,r}(x) \right],$$
for $x\in\mathbb{R}$ such that the series
\Eq{E}
$$E^k_{j,r}(x)=\sum_{n=0}^\infty \frac{e^{i(nk+j)x}}{(nk+j)^r}$$
is convergent. 
We call \emph{nodes} the points $x=2\pi l/k$ for $l\in\mathbb{Z}$, where the series diverges.
   
The reason for the terminology ``trigonometric polylogarithm'' is that, as we discuss in the Appendix,  these functions are finite linear combinations of polylogarithms evaluated at trigonometric arguments.
It is readily seen that, for $r\geq 1$,
\begin{equation} \label{eqn:SCdefn}
    S^k_{j,r}(x) = \sum_{n=0}^\infty \frac{\sin[(nk+j)x]}{(nk+j)^r}, \qquad 
    C^k_{j,r}(x) = \sum_{n=0}^\infty \frac{\cos[(nk+j)x]}{(nk+j)^r}.
\end{equation} 

The main properties of these functions are summarised in the following proposition, which is proved in Appendix~\ref{sec:ProofOfPropSC}.

\begin{proposition}\label{propSC}
    For all $r\geq 1$, at all points other than the nodes, the trigonometric polylogarithms $S^k_{j,r}$ and $C^k_{j,r}$ are $\ro C^\infty$ smooth and are arranged in derivative chains
        \begin{equation} \label{eqn:SCDerivatives}
        \frac{\d}{\d x}C^k_{j,r+1}(x) = -S^k_{j,r}(x), \qquad \frac{\d}{\d x}S^k_{j,r+1}(x) = (-1)^{r}C^k_{j,r}(x).
    \end{equation}
    For $r\geq2$, the trigonometric polylogarithms $S^k_{j,r}$ and $C^k_{j,r}$ defined above belong to $\ro C_{\mathrm{per}}^{r-2}[-\pi,\pi]$, the space of continuous periodic functions with $r-2$ continuous periodic derivatives.
    \hfill\break \sstrut1 \qquad 
    
    \noindent
    {\small \bf Behaviour at the nodes}
    
\smallskip
   \noindent For $r=1$, at the nodes $x=2\pi l/k$, for $l\in\mathbb{Z}$,
    \begin{multline*}
        \left.
        \begin{matrix}S^k_{j,1}(x)\sstrut9\\C^k_{j,1}(x)\end{matrix}
        \right\}
        \mbox{ has a jump of height }
        \left\{\begin{matrix}\dsty\frac{\pi}{k}\cos\left(\frac{2\pi jl}{k}\right)\sstrut{15}\\\dsty\frac{\pi}{k}\sin\left(\frac{2\pi jl}{k}\right)\end{matrix}\right\}\dsty
        \mbox{ superimposed with a logarithmic}
        \\
        \mbox{cusp pointing up/down/not present if }
        \left\{
        \begin{matrix}\dsty\sin\left(\frac{2\pi jl}{k}\right)\sstrut{15}\\\dsty\cos\left(\frac{2\pi jl}{k}\right)\end{matrix}
        \right\}
        \mbox{ is positive/negative/zero.}
    \end{multline*}
    For $r=2$, at the nodes $x=2\pi l/k$, for $l\in\mathbb{Z}$,
    \begin{multline*}
        \left.
        \begin{matrix}S^k_{j,2}(x)\sstrut9\\C^k_{j,2}(x)\end{matrix}
        \right\}
        \mbox{ has a point of infinite gradient,}
        \mbox{ unless }
        \left\{
        \begin{matrix}\dsty\cos\left(\frac{2\pi jl}{k}\right)=0\sstrut{15}\\\dsty\sin\left(\frac{2\pi jl}{k}\right)=0\end{matrix}
        \right\}\\
        \mbox{ in which case there is a corner.}\hfill
    \end{multline*}
\end{proposition}

The most important case of trigonometric polylogarithmic functions for the present work is the case $r=1$, and we will at times refer to these as {\em trigonometric hypergeometric functions}; see the Appendix. In this case we use the notation
\Eq{Sjk1}
$$S^k_j(x) = S^k_{j,1}(x) = \Sumoi n\frac{\sin(nk+j)x}{nk+j}.$$
In the Appendix, we provide the alternative formula
\begin{multline} \label{Skjl}
    S^k_j(x) = \frac{1}{k}\sum_{l=1}^k \left[ \sin\left(\frac{2\pi j\:l}{k}\right)\log\left\lvert2\sin\left(\frac{x}{2}+\frac{\pi l}{k}\right)\right\rvert \right. \\ \left. + \cos\left(\frac{2\pi j\:l}{k}\right)\frac{\sign\left(x+\fraco{2\pi l}{k}\right)\pi-\left(x+\fraco{2\pi l}{k}\right)}{2} \right].
\end{multline}

\begin{figure}[h]
  \begin{subfigure}{.3\textwidth}
\centering
       \includegraphics[height=1.3in]{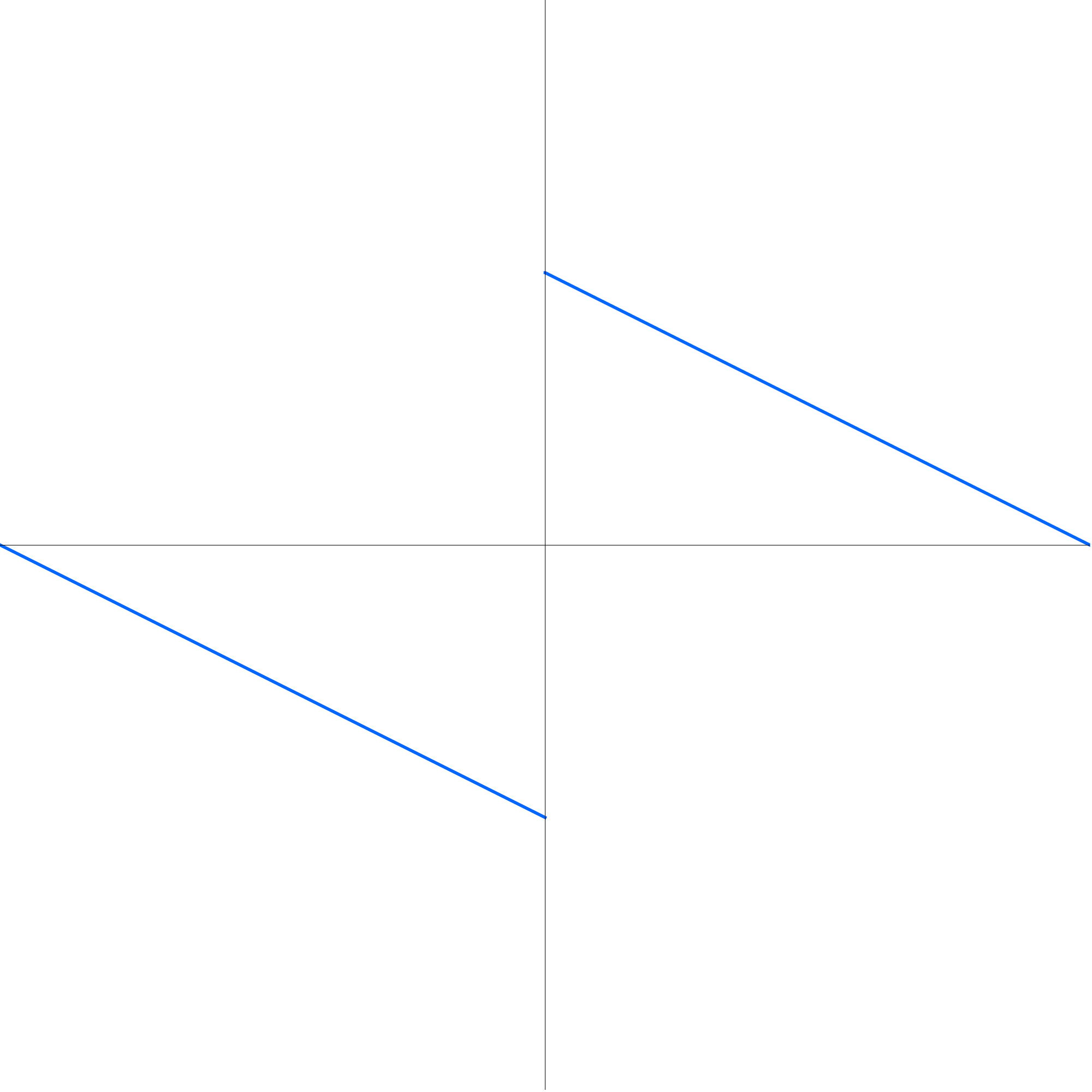}
       \subcaption*{$S^1_1(x)$}
           \end{subfigure}
            \begin{subfigure}{.3\textwidth}
\centering
       \includegraphics[height=1.3in]{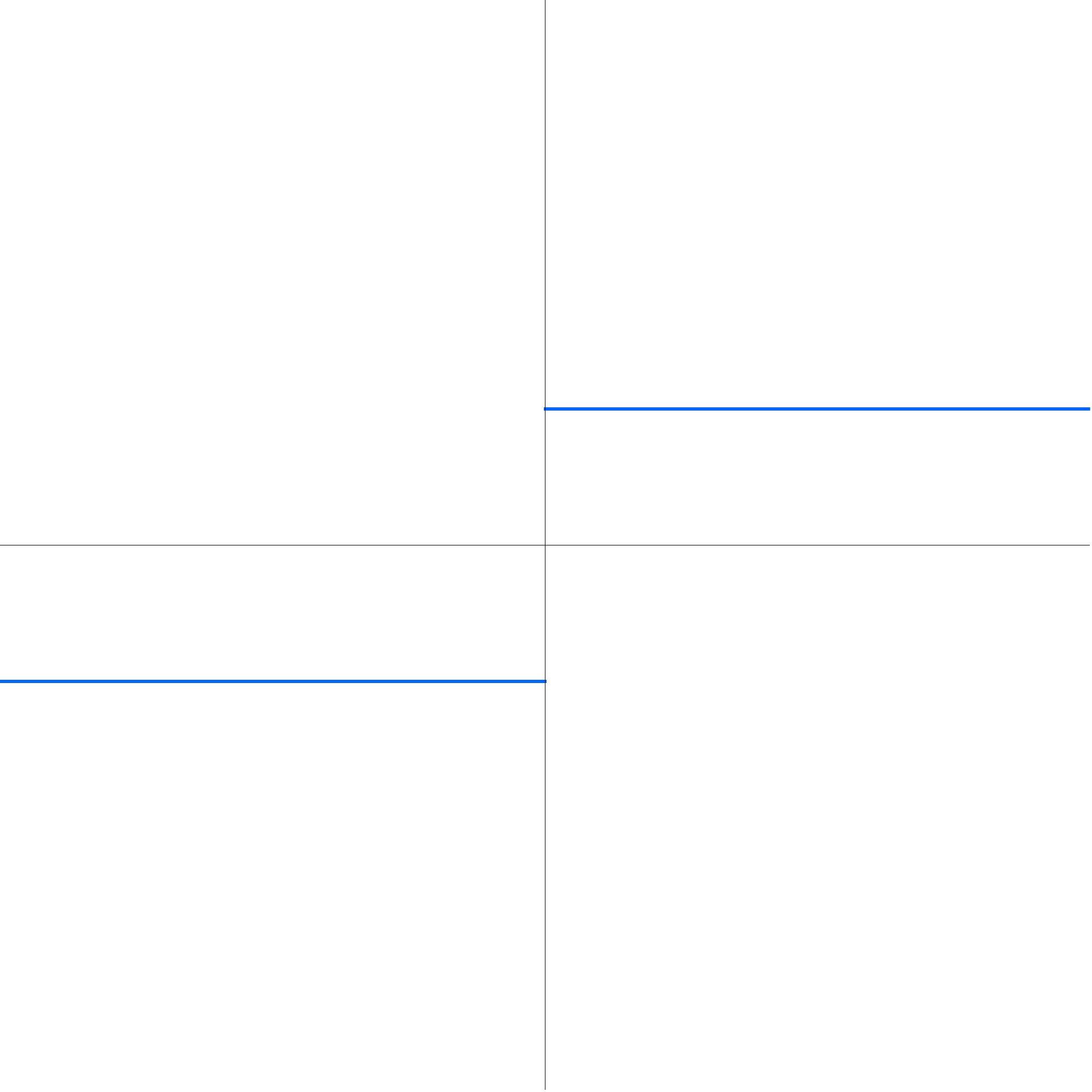}
    \subcaption*{$S^2_1(x)$}
           \end{subfigure}
              \begin{subfigure}{.3\textwidth}
\centering
       \includegraphics[height=1.3in]{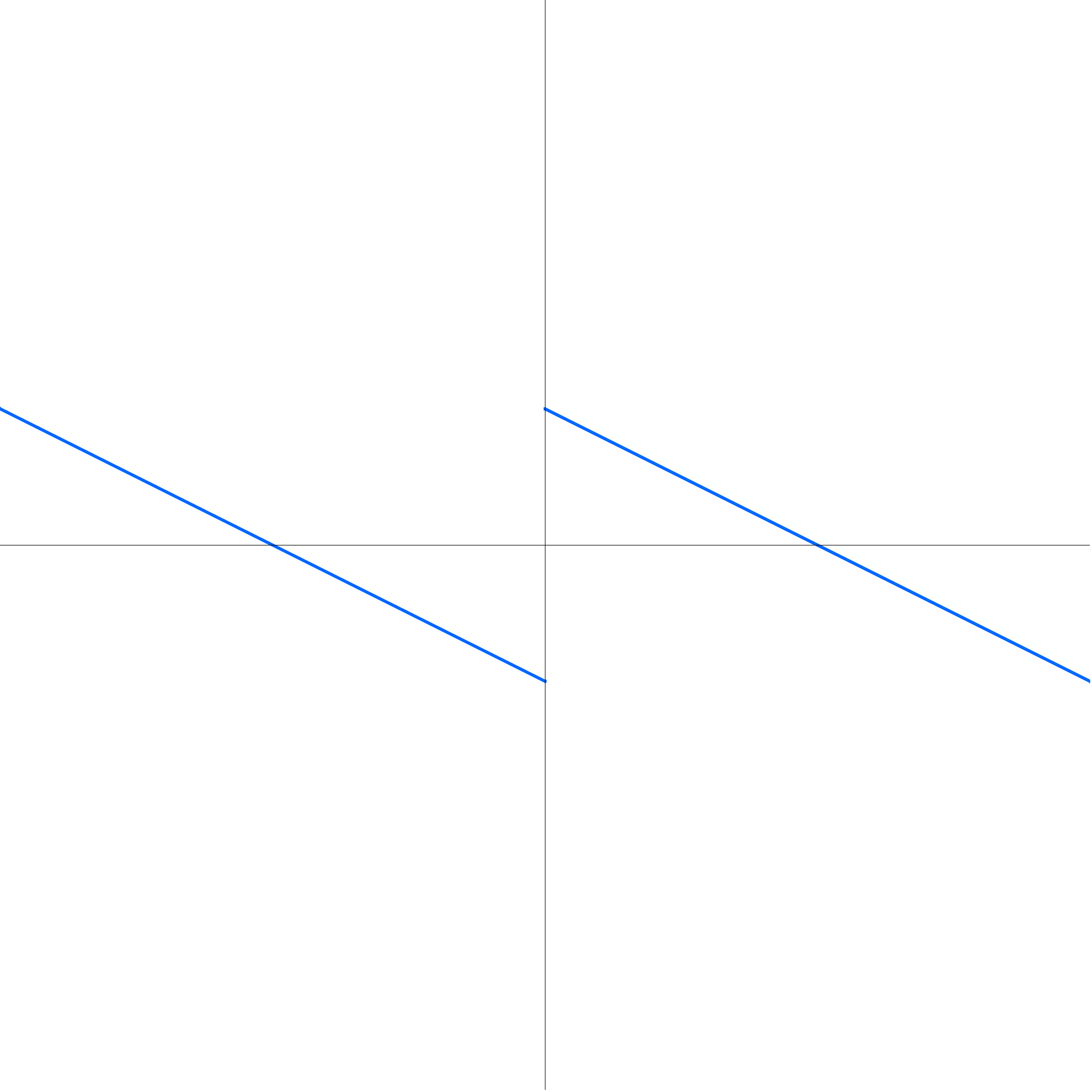}
    \subcaption*{$S^2_2(x)$}
           \end{subfigure}
              \begin{subfigure}{.3\textwidth}
\centering
       \includegraphics[height=1.3in]{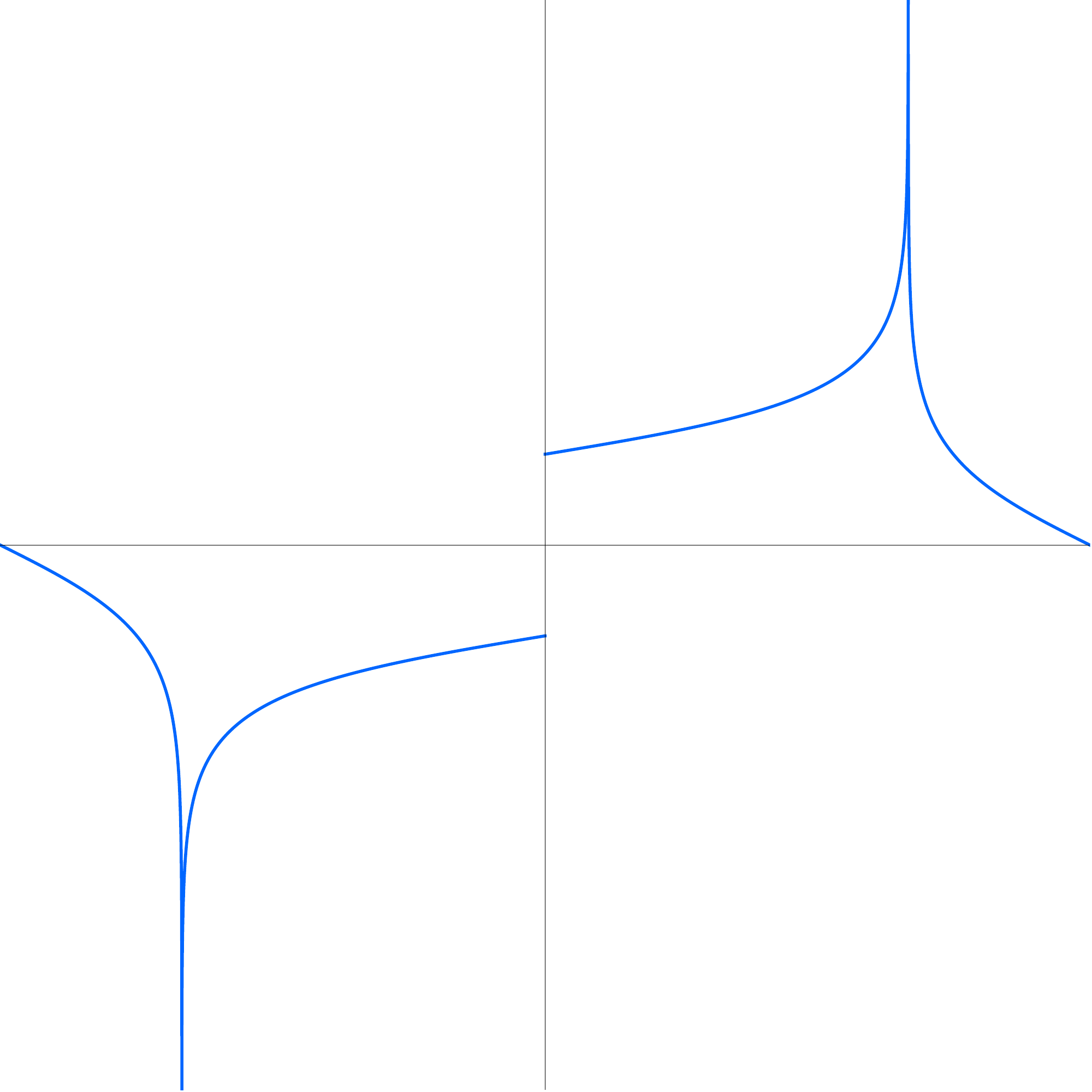}
    \subcaption*{$S^3_1(x)$}
           \end{subfigure}\hglue.28in
               \begin{subfigure}{.3\textwidth}
\centering
       \includegraphics[height=1.3in]{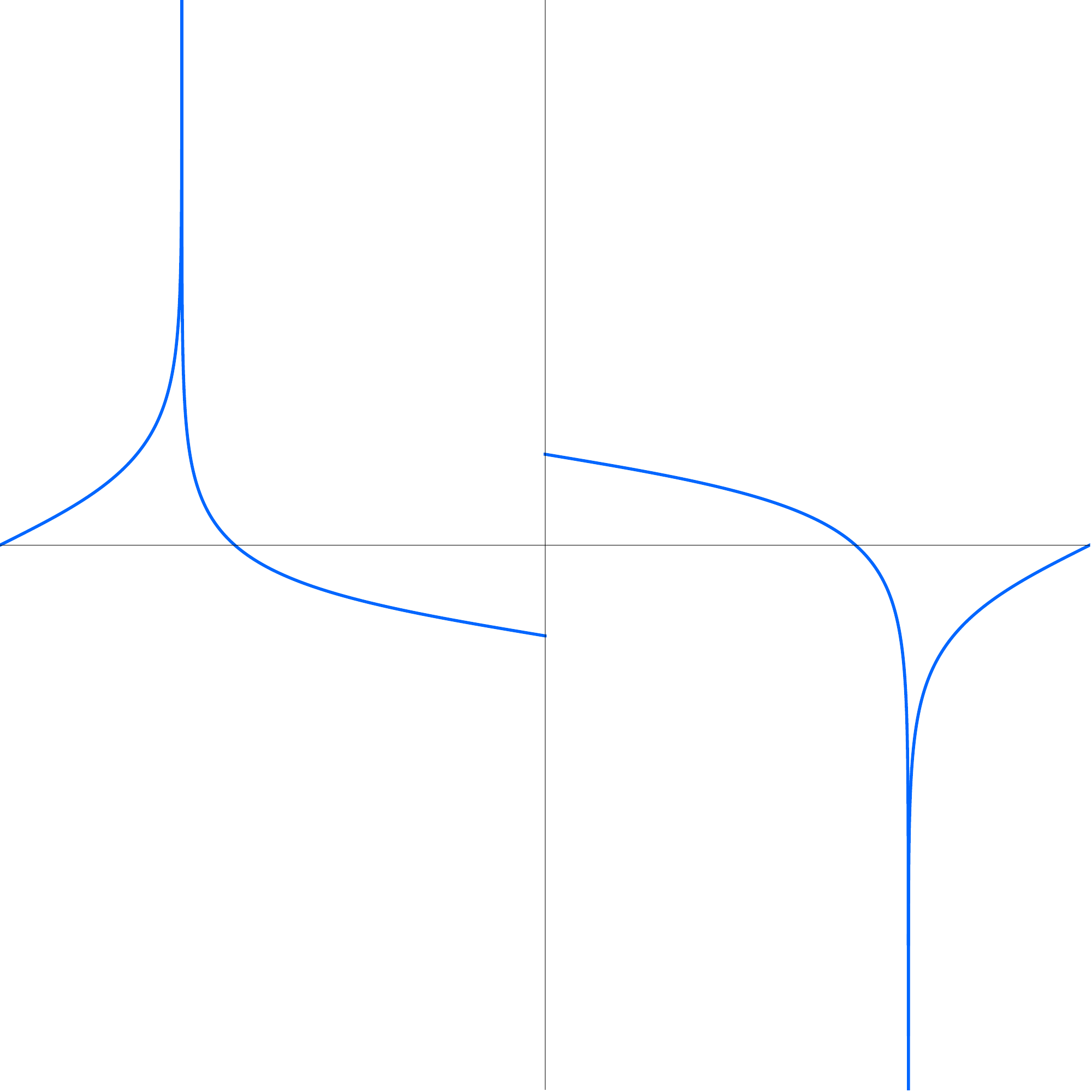}
    \subcaption*{$S^3_2(x)$}
           \end{subfigure}\hglue.28in
               \begin{subfigure}{.3\textwidth}
\centering
       \includegraphics[height=1.3in]{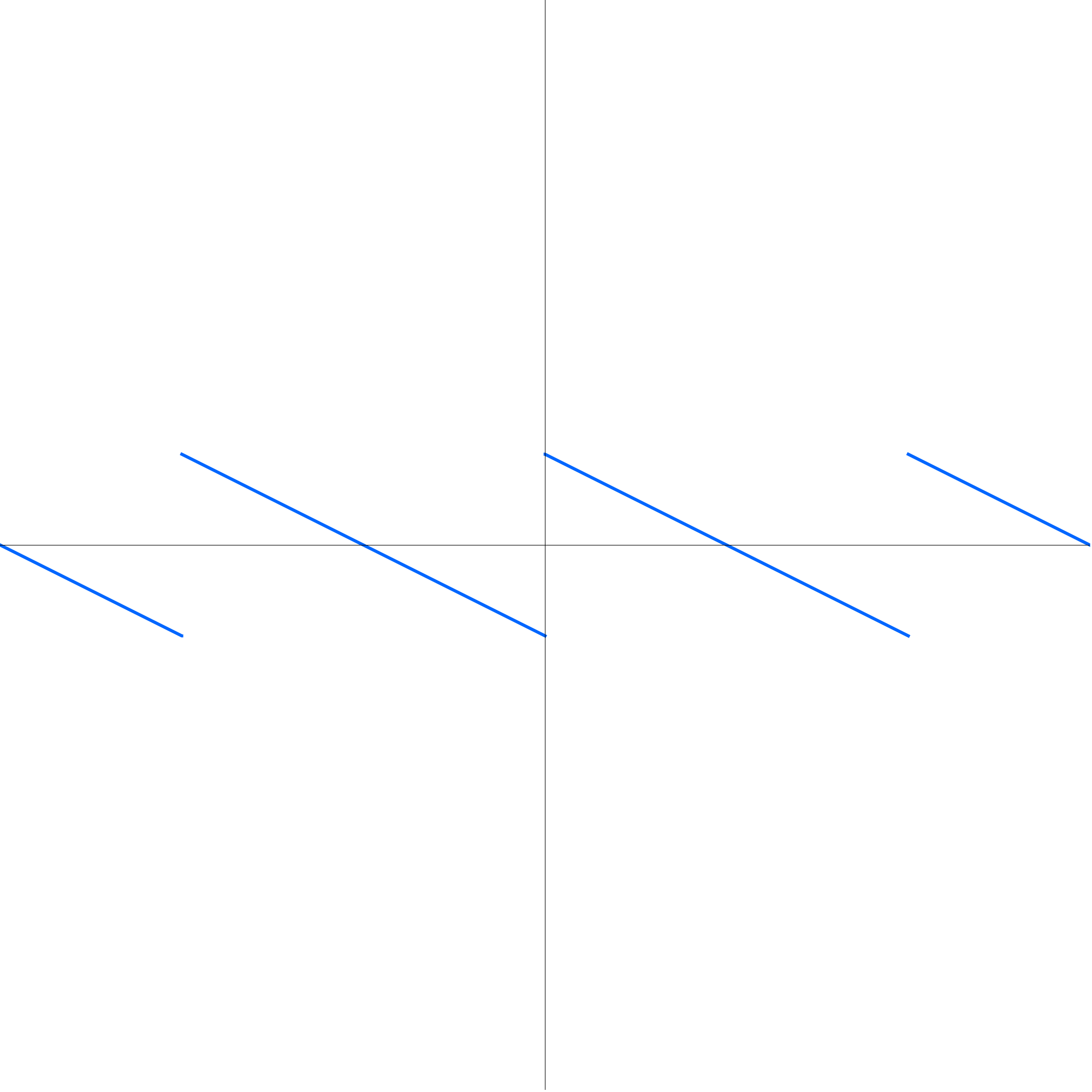}
    \subcaption*{$S^3_3(x)$}
           \end{subfigure}
        \caption{Graphs of trigonometric hypergeometric functions $S^k_j(x)$, $k=1,2,3$. \\ \hglue1in 
            Both  the horizontal and vertical axes are from $-\pii$ to $\pii$.}
        \label{fig:Skj.plots2.1}
            \end{figure}

\begin{figure}[t]
            \begin{subfigure}{.3\textwidth}
\centering
       \includegraphics[height=1.3in]{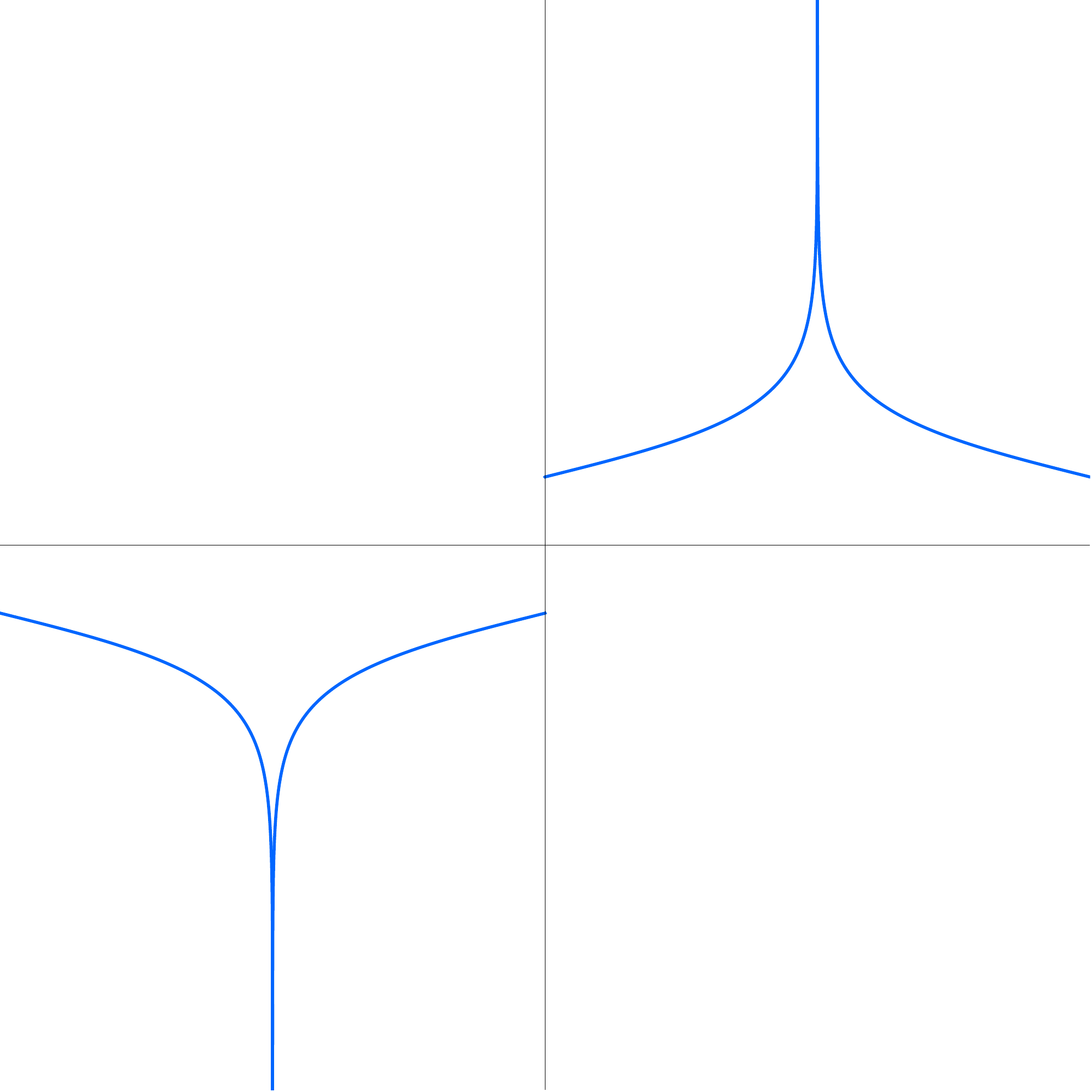}
    \subcaption*{$S^4_1(x)$}
           \end{subfigure}
               \begin{subfigure}{.3\textwidth}
\centering
       \includegraphics[height=1.3in]{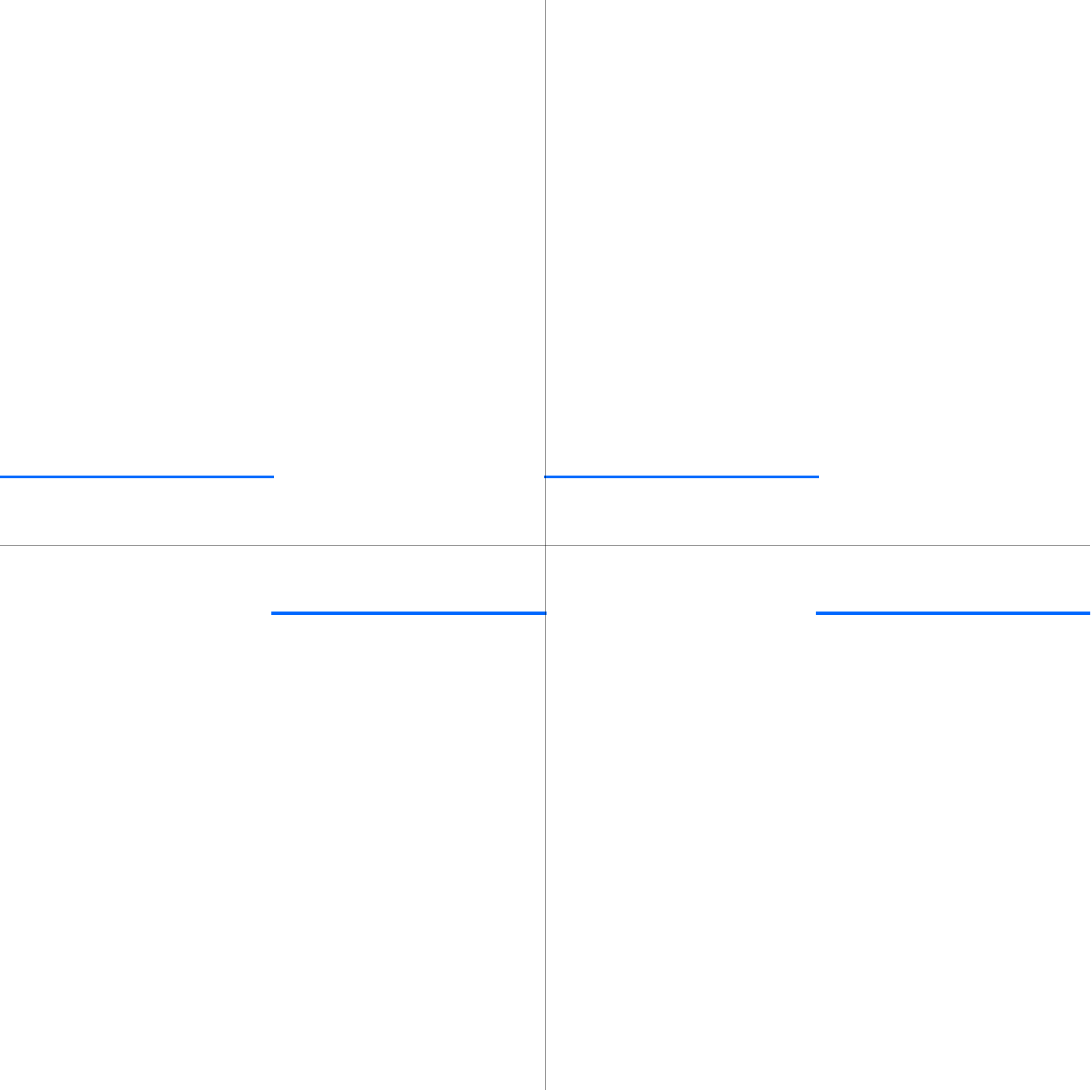}
    \subcaption*{$S^4_2(x)$}
           \end{subfigure}
               \begin{subfigure}{.3\textwidth}
\centering
       \includegraphics[height=1.3in]{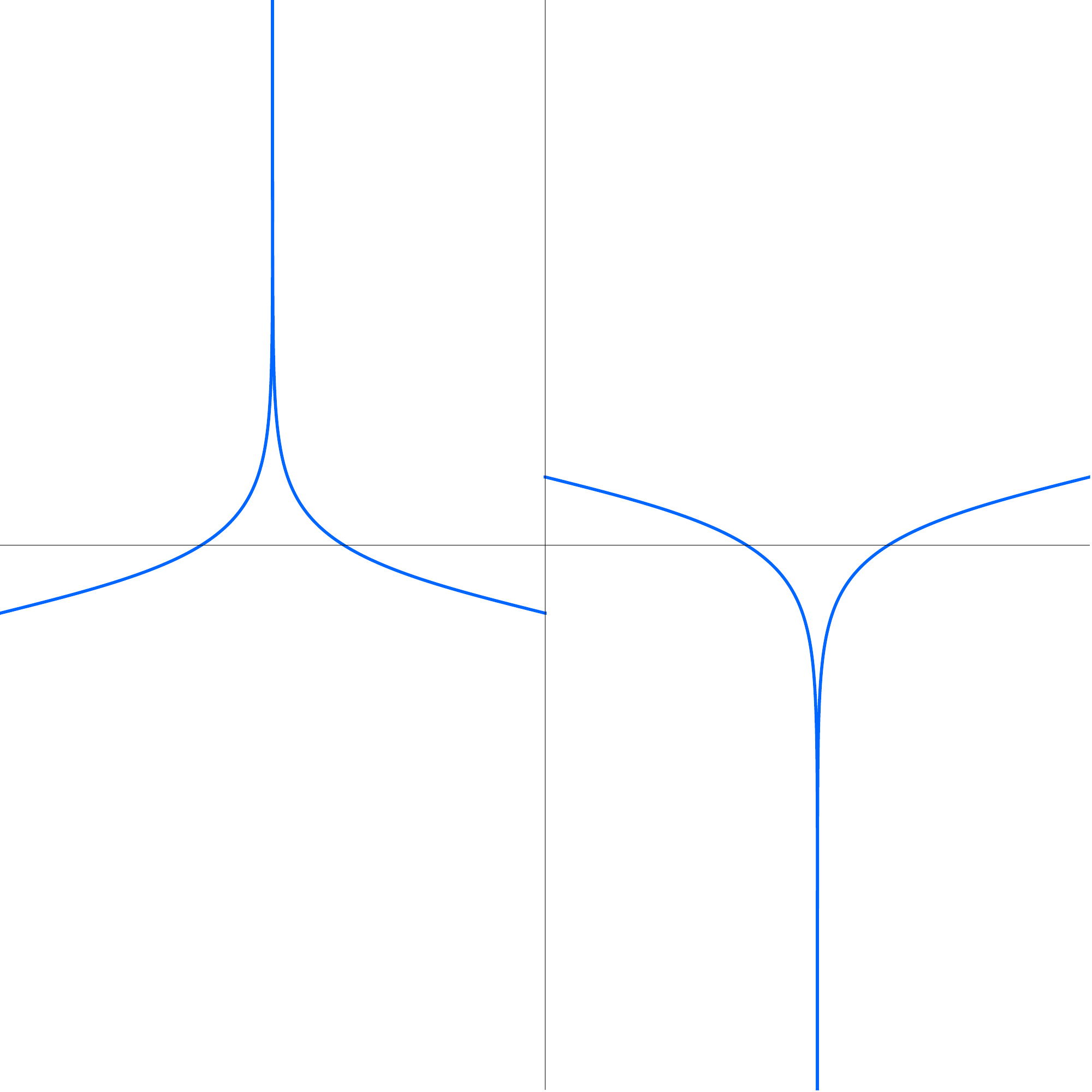}
    \subcaption*{$S^4_3(x)$}
           \end{subfigure}
               \begin{subfigure}{.3\textwidth}
\centering
       \includegraphics[height=1.3in]{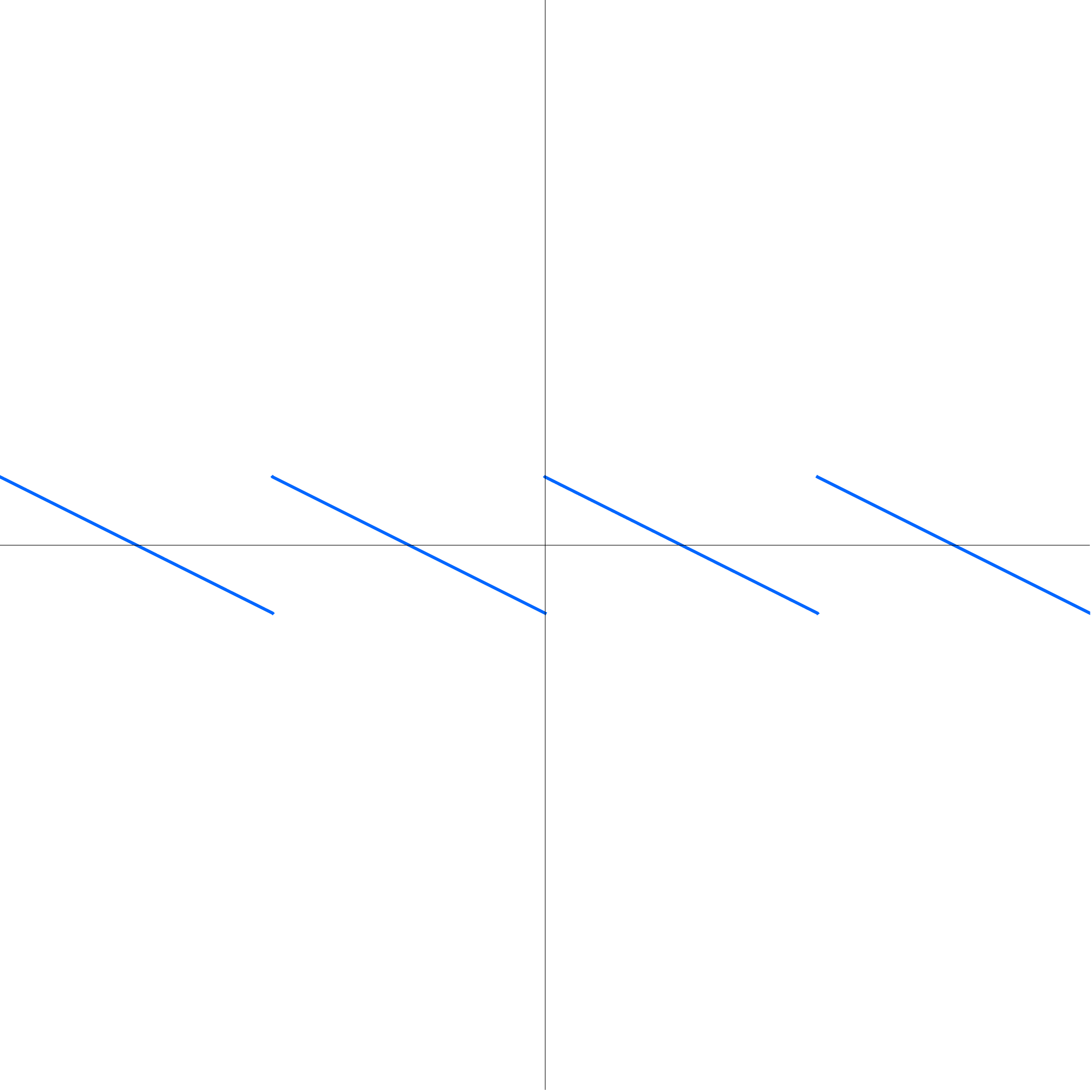}
    \subcaption*{$S^4_4(x)$}
           \end{subfigure}
               \begin{subfigure}{.3\textwidth}
\centering
       \includegraphics[height=1.3in]{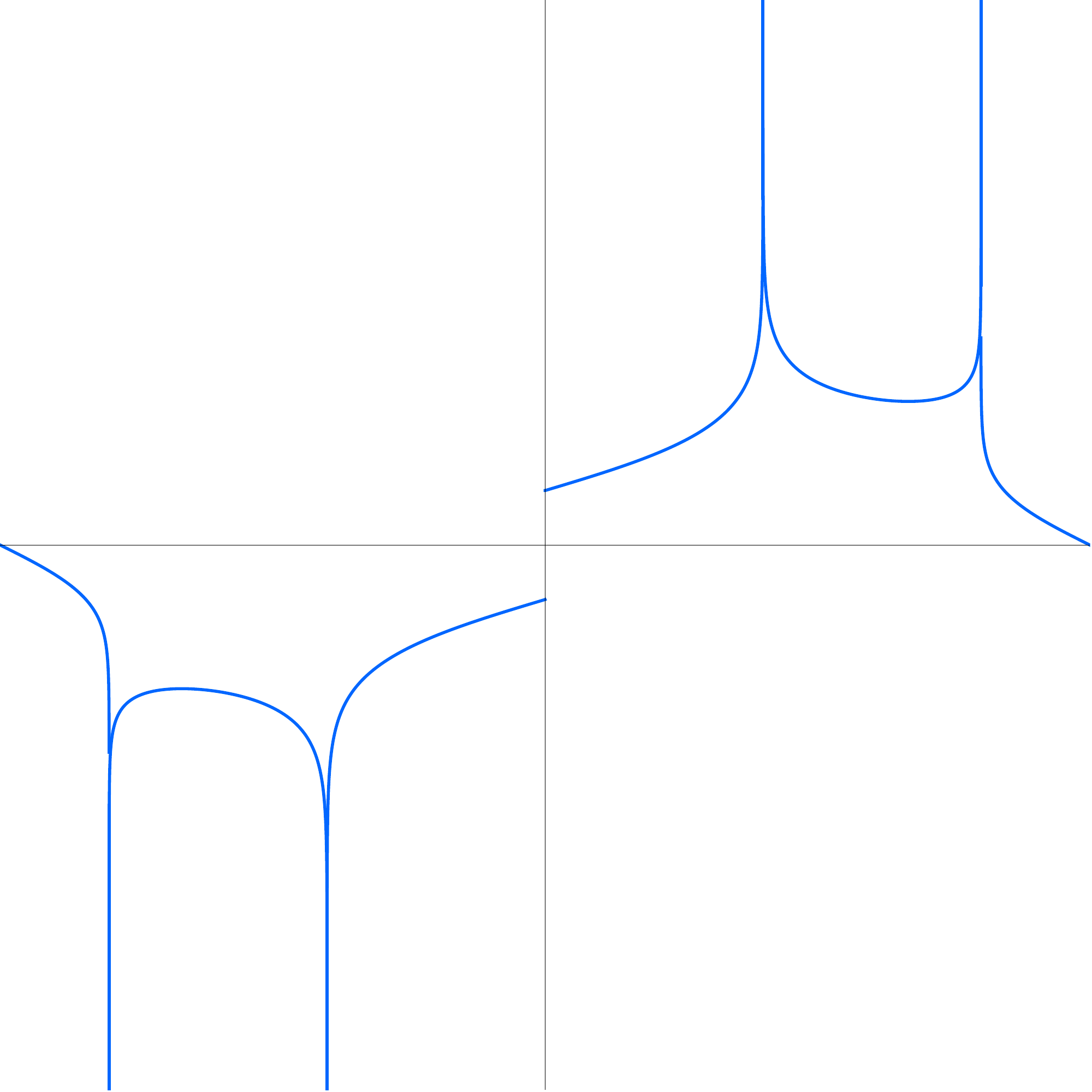}
    \subcaption*{$S^5_1(x)$}
           \end{subfigure}
               \begin{subfigure}{.3\textwidth}
\centering
       \includegraphics[height=1.3in]{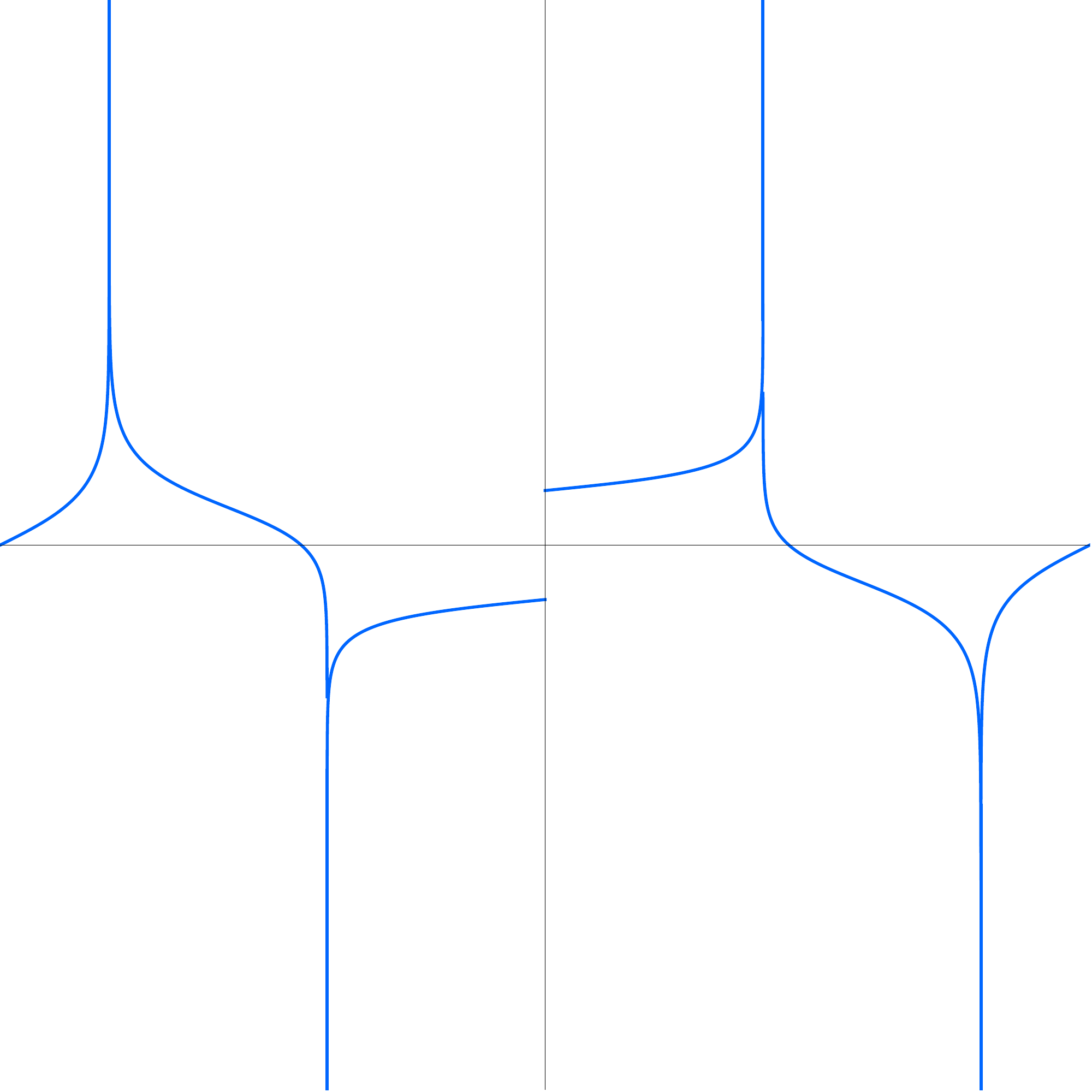}
    \subcaption*{$S^5_2(x)$}
           \end{subfigure}
               \begin{subfigure}{.3\textwidth}
\centering
       \includegraphics[height=1.3in]{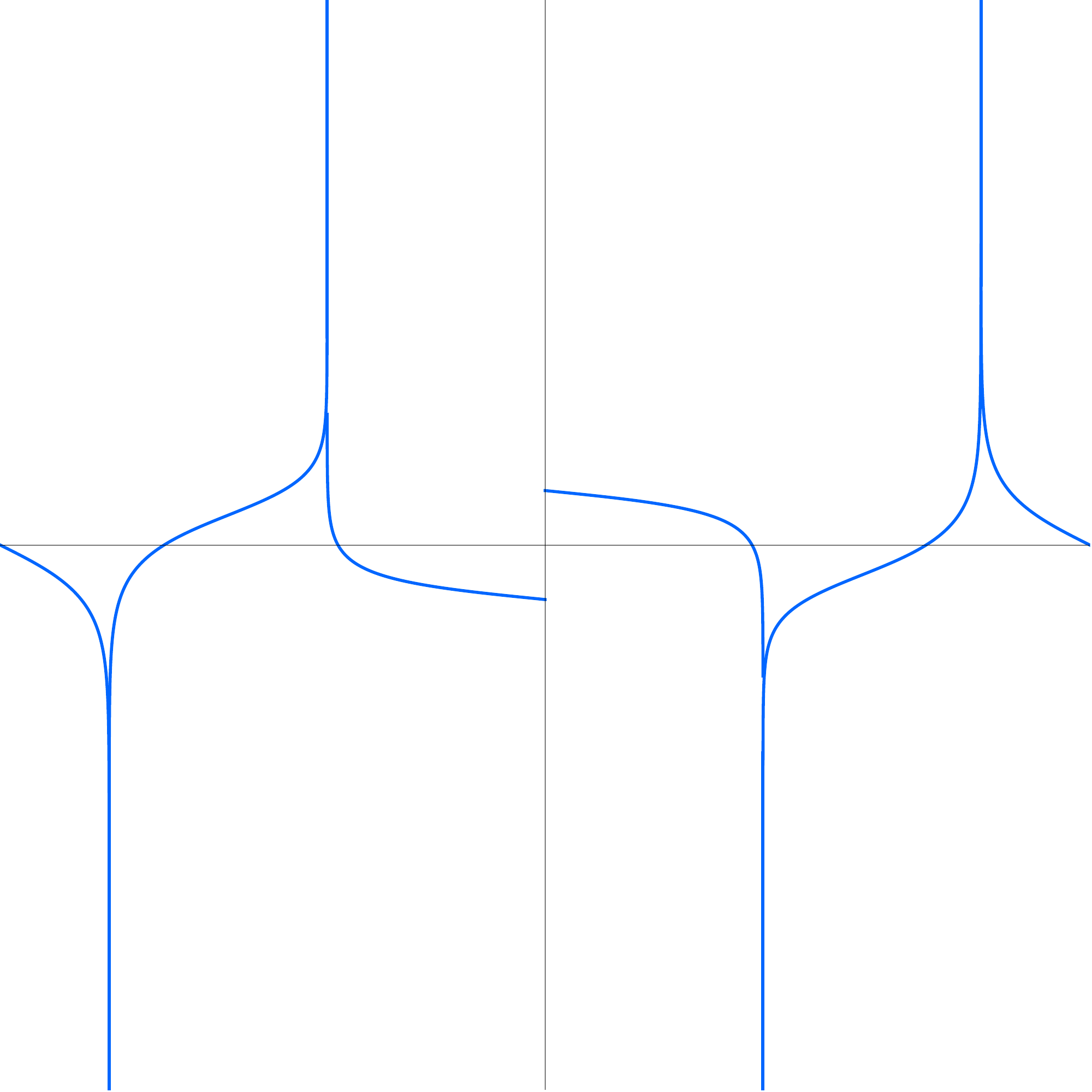}
    \subcaption*{$S^5_3(x)$}
           \end{subfigure}\hglue.3in
               \begin{subfigure}{.3\textwidth}
\centering
       \includegraphics[height=1.3in]{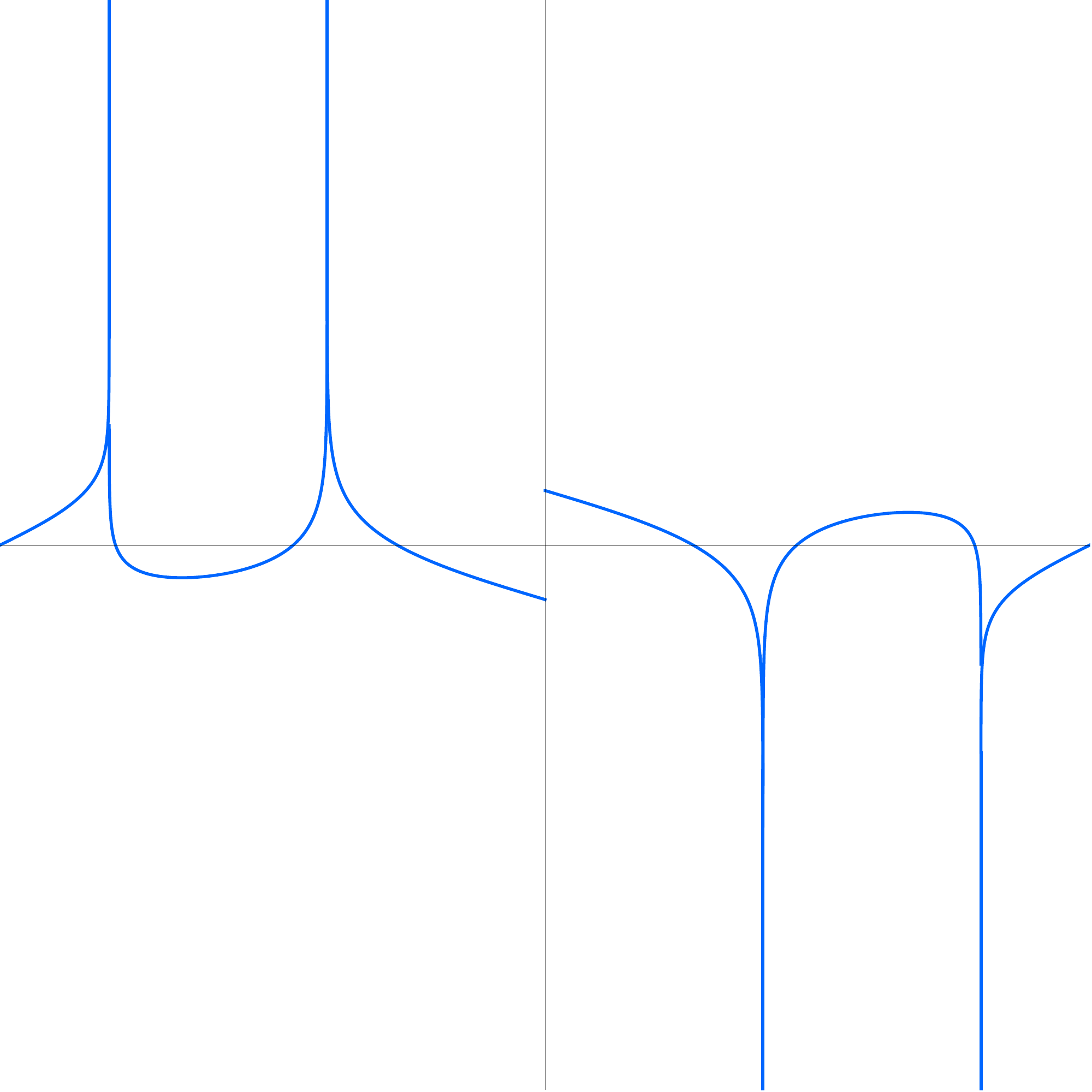}
    \subcaption*{$S^5_4(x)$}
           \end{subfigure}\hglue.3in
               \begin{subfigure}{.3\textwidth}
\centering
       \includegraphics[height=1.3in]{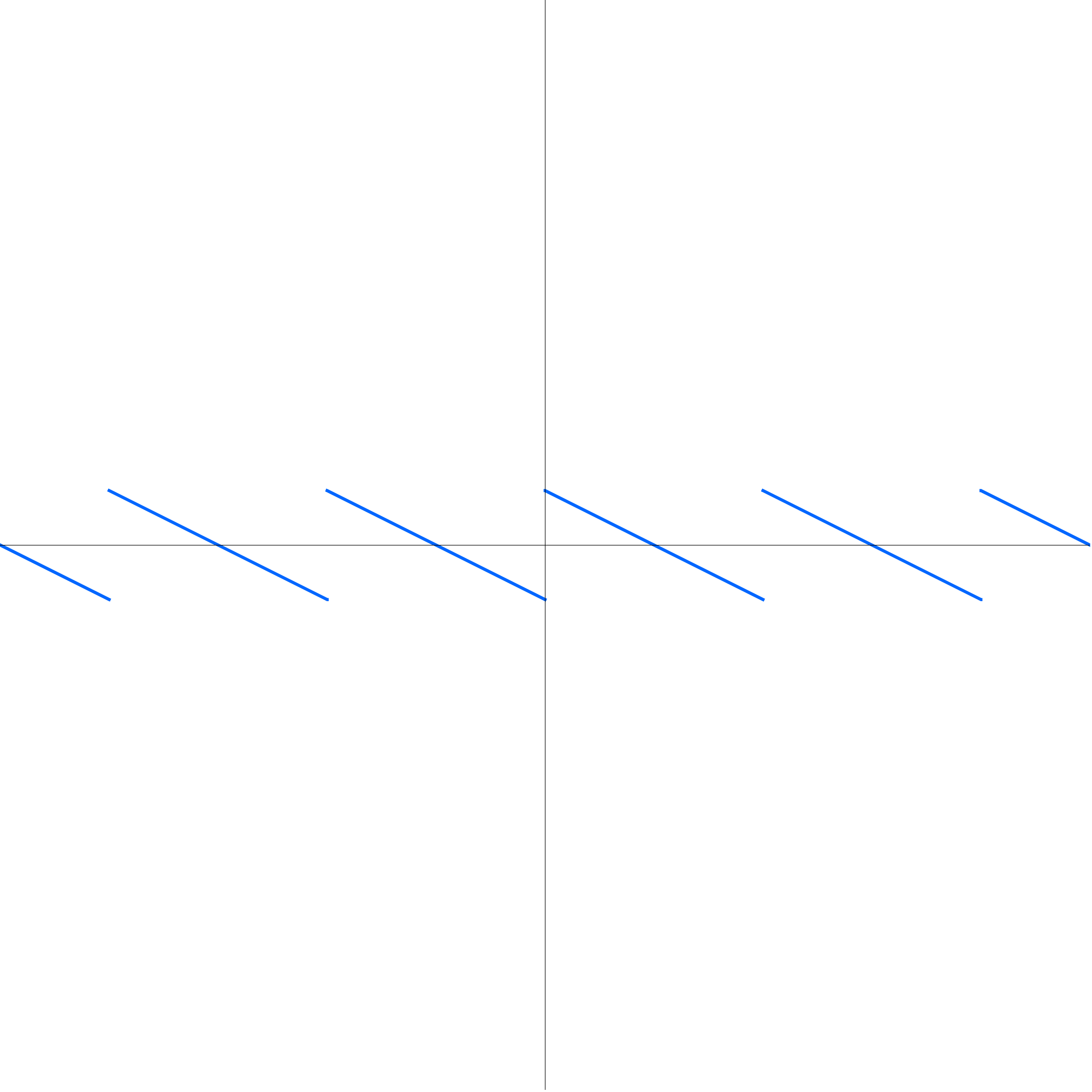}
    \subcaption*{$S^5_5(x)$}
           \end{subfigure}
        \caption{Graphs of trigonometric hypergeometric functions $S^k_j(x)$, $k=4,5$. \\ 
        \hglue1in 
       Both  the horizontal and vertical axes are from $-\pii$ to $\pii$.}
        \label{fig:Skj.plots2.2}
\end{figure}

The discontinuities of these functions are described in \pr{propSC} and the low order ones are displayed in Figure~\ref{fig:Skj.plots2.1} and \ref{fig:Skj.plots2.2}.  
In the graphs, which were plotted in \Mathematica, the infinite logarithmic cusps have been drawn in by hand, as the plotting routines made them misleadingly appear to be finite in height.


These functions have distributional derivatives
\Eq{dSkjl}
$$\eeq{\od{S^k_j}x  =
\frac\pi k\,\Sumo l{k-1}\cos\Paz{\frac{2\pii \:l\:j}k} \dide \Paz{x+ \frac{2\pii \:l}k} +  \fra{2\:k} \Sum l{k-1} \sin\Paz{\frac{2\pii \:l\:j}k}\cot\Paz{\f2\:x + \frac{\pii \:l}k} .}$$
The first summation in expression \eq{dSkjl} is a linear combination of Dirac deltas at the nodes and the second is a linear combination of pole singularities at the nodes that are not integer multiples of $\pi$.
Furthermore, for $\rg j{k-1}$, the sum
\Eq{Skjsum}
$$S^k_j(x) + S^k_{k-j}(x) = \fra{2\i}\Sumii n\Pa{\frac{e^{\i(nk+j)x}}{nk+j} + \frac{e^{\i(nk-j)x}}{nk-j}}$$
is a piecewise constant function on the intervals $m\pii/k < x < (m+1)\pii/k$; the proof relies on the method used in \rf{Odq} to characterize Fourier series representing piecewise constant functions. 
Note that the pole terms cancel out in the sum
\Eq{dSkjsum}
$$\od{S^k_j}x + \od{S^k_{k-j}}x = \frac{2\pii} k\,\Sumo l{k-1}\cos\Paz{\frac{2\pii \:l\:j}k} \dide \Paz{x+ \frac{2\pii \:l}k}$$
leaving a linear combination of deltas at the nodes, which is in accordance with our earlier observation that the sum \eq{Skjsum}
 is piecewise constant.  Indeed, \eq{dSkjsum} prescribes the magnitudes of its jump discontinuities, each given by the coefficient of the corresponding Dirac deltas.

\section{New Revival Phenomena}
\label{r}

In all that follows we consider the solution of the \emph{periodic Riemann problem}, posed on the interval $(-\pi,\pi)$,  for each of the linear equations of the section~\ref{sec:ModelWaveEquations}. In other words, we study the periodic initial-boundary value problem problem obtained when the initial condition is the unit step function
\be
u(0,x) = \frac{1 + \sign x}2,\qquad x\in(-\pi,\pi).
\label{riemannIC}
\ee
Without further mention, here and elsewhere we assume that functions and distributions defined in $(-\pi,\pi)$ are extended $2\pi$-periodically to $\mathbb{R}$ in the usual way, if required.

\smallskip
Assuming that the dispersion relation $\omega (k)$ is an odd function\footnote{One can easily write down the solution for even and more general dispersion relations, but the formula is more complicated, and all examples considered here are odd.} of $k$, the solution to the initial value problem, for a linear integro-differential equation  of the form \eq{leq}, \eq{leqF} is formally given by, \crf{COdisp},
\Eq{sol}
$$u(t,x) = \fra 2 + \frac 2\pi \Sumoi j \frac{\sin\bk{(2\:j+1)x - \omega (2\:j+1)t}}{2\:j+1}.$$

Our main result  is that the solution of the Riemann problem for these equations, at times equal to rational multiples of the period, is given by a linear superposition of translated trigonometric  hypergeometric functions or of translated trigonometric polylogarithms. For BO, this is an exact result, summarised in \th{lbof} below. 

\begin{rem}
All figures in this paper are generated by \Mathematica\ and plotted with horizontal axis $[-\pi,\pi]$ and vertical axis $[-2,2]$.  The cusped behavior at rational times  that is shown on the figures is consistent with our analysis, although \Mathematica\ is unable to plot the logarithmic cusps at infinite height, misleadingly indicating that they are finite in extent.  Unlike Figures~\ref{fig:Skj.plots2.1} and \ref{fig:Skj.plots2.2}, we have chosen not to draw the infinitely extended, but vanishingly narrow cusps here, leaving it to the reader to fill in what the actual profile should look like.
\end{rem}
\subsection{Revivals for the Linearised Korteweg--deVries Equation} \label{ssec:lKdV}

Before investigating the equations that form the focus of this work, let us review the known revival phenomena,  \cite{Odq}, for the periodic linearised Korteweg--deVries (KdV) equation
\be
u_t+u_{xxx}=0.
\label{lkdv}\ee
This linear equation is also known  as {\em Airy equation} since its fundamental solution on the real line can be expressed in terms of the Airy function.

The dispersion relation  is $\omega_{kdv}(k) = k^3$, which is a polynomial with integer coefficients.  Thus, the general result of \cite {Odq} (see also \rf{COdisp, ET}) yields the following: 

\Th{lkdvf} At a rational time $2\pi\:p/q$, the solution to the periodic initial-boundary value problem for the linearised KdV equation (\ref{lkdv}) on the interval $-\pii < x < \pi$, with initial datum $u(0,x) = u_0(x)$, is a linear combination of translates $u_0(x - 2\pi\:j/q)$, $\rgo j{q-1}$:
$$
u(2\pi\:p/q,x)=\sum_{j=0}^{q-1}a_{p,q}(j)u_0(x-2\pi \:j/q)
$$
for coefficients $a_{p,q}(j)\in \mathbb{R}$ explicitly computable.

\vspace{.15in}
\begin{figure}[h]
  \begin{subfigure}{.3\textwidth}
\centering
       \includegraphics[height=.8in]{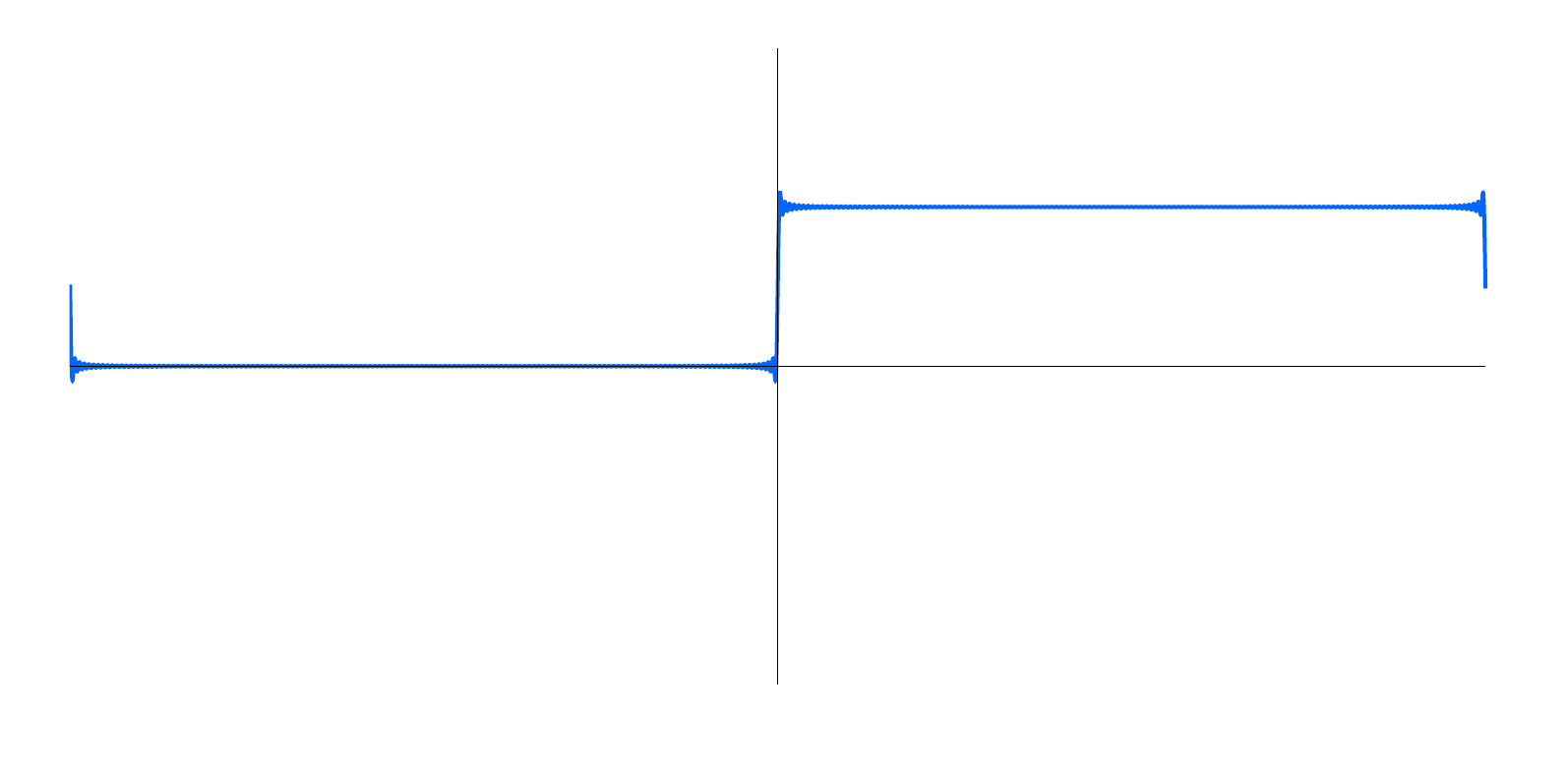}
       \subcaption*{$t=0$}
           \end{subfigure}
            \begin{subfigure}{.3\textwidth}
\centering
       \includegraphics[height=.8in]{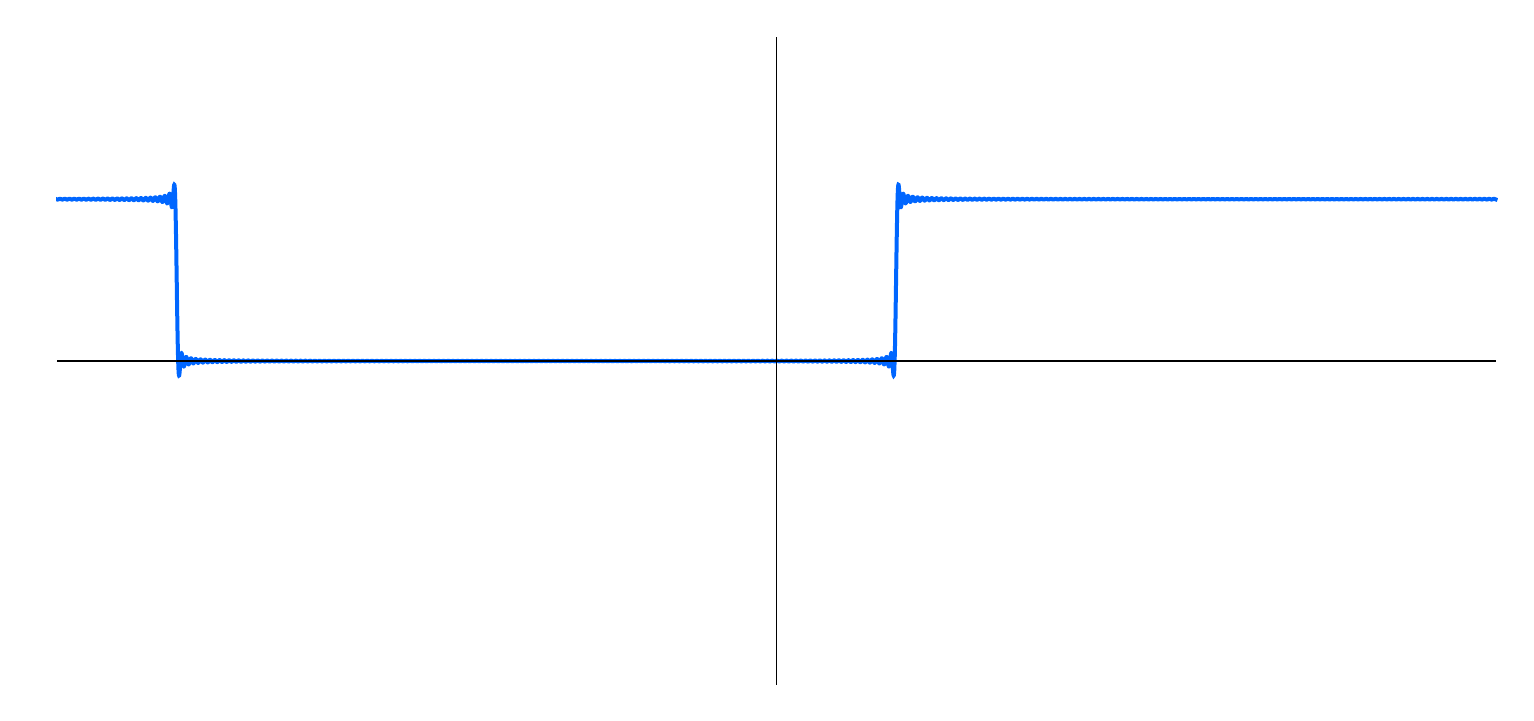}
    \subcaption*{$t=\frac \pi 6$}
           \end{subfigure}
              \begin{subfigure}{.3\textwidth}
\centering
       \includegraphics[height=.8in]{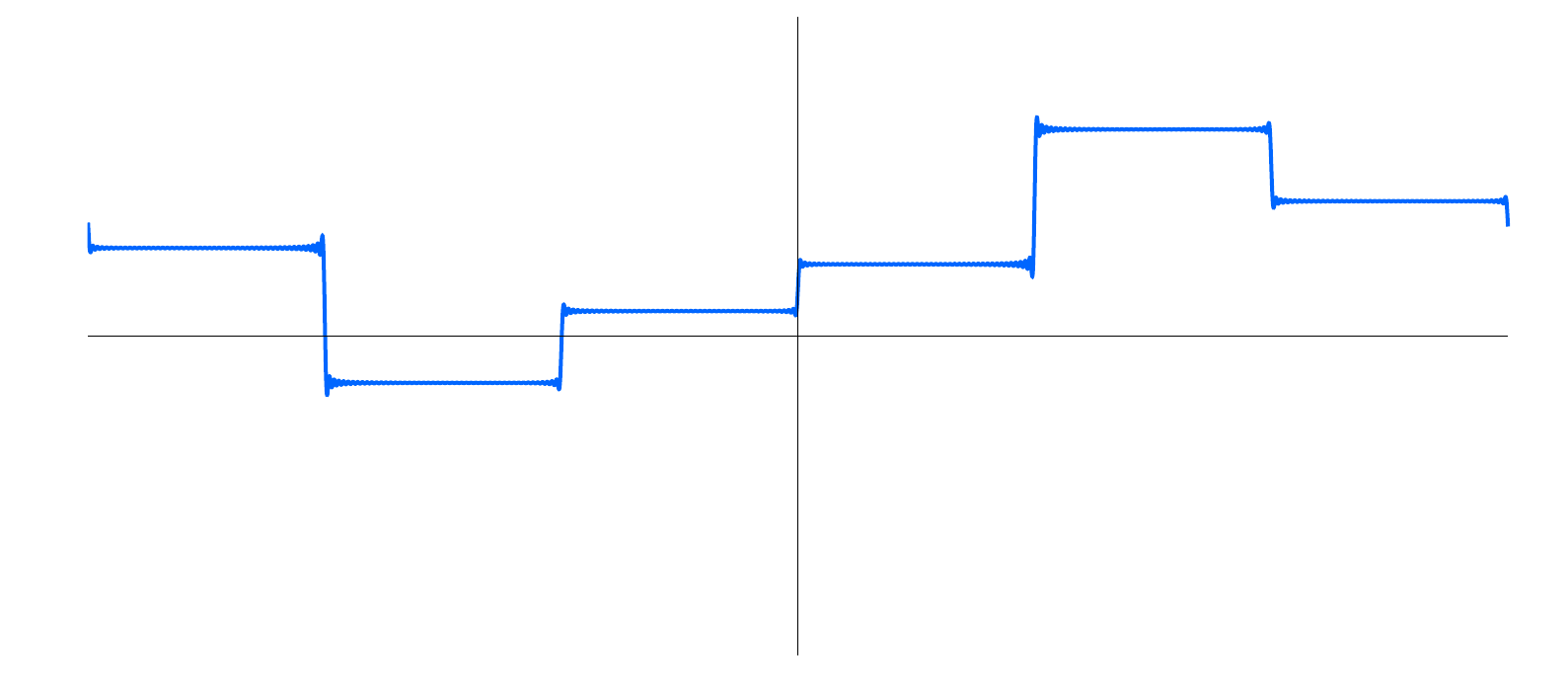}
    \subcaption*{$t=\frac \pi 9$}
           \end{subfigure}
 \caption{The solution of the periodic Riemann problem for the linearised KdV equation at some representative times that are rational multiples of $\pi$ -  Note the effect of Gibbs phenomenon at points of discontinuity. The horizontal  axis is from $-\pii$ to $\pii$, the vertical axis is from $-2$ to $2$} \label{kdvfig}
\end{figure}

In particular, the solution to the periodic Riemann problem given in the general form \eq{sol} as
\Eq{kdvstep}
$$u(t,x) = \fra 2 + \frac 2\pi \Sumoi j \frac{\sin\bbk{(2\:j+1)\bpa{x - (2\:j+1)^2t}}}{2\:j+1}.$$
 evaluated at rational times $t = \pi p/q$ is a linear combination of translated step functions, and hence piecewise constant, being constant on the intervals $\ivl{\pi\:j/q}{\pi\:(j+1)/q}$. The solution may be constant on longer subintervals; see \rf{OT} for a number-theoretic characterisation of these occurrences.

On the other hand, at irrational times, the solution profile is a continuous fractal, with fractal dimension $D$ bounded by $\fr54 \leq D \leq \fr74$; see \rf{OskolkovV,ErSh}.

\subsection{Revival for the Linearised Benjamin--Ono Equation}
\label{rBO}

In this section we make use of the properties of  trigonometric hypergeometric functions  to study revival properties in the case of the linearised Benjamin--Ono equation \eq{lBO} subject to periodic boundary conditions on the interval $\ipp$.

\begin{figure}[h]
  \begin{subfigure}{.5\textwidth}
\centering
       \includegraphics[height=1.2in]{./BOpi0}
       \subcaption*{$t=0$}
           \end{subfigure}
            \begin{subfigure}{.5\textwidth}
\centering
       \includegraphics[height=1.2in]{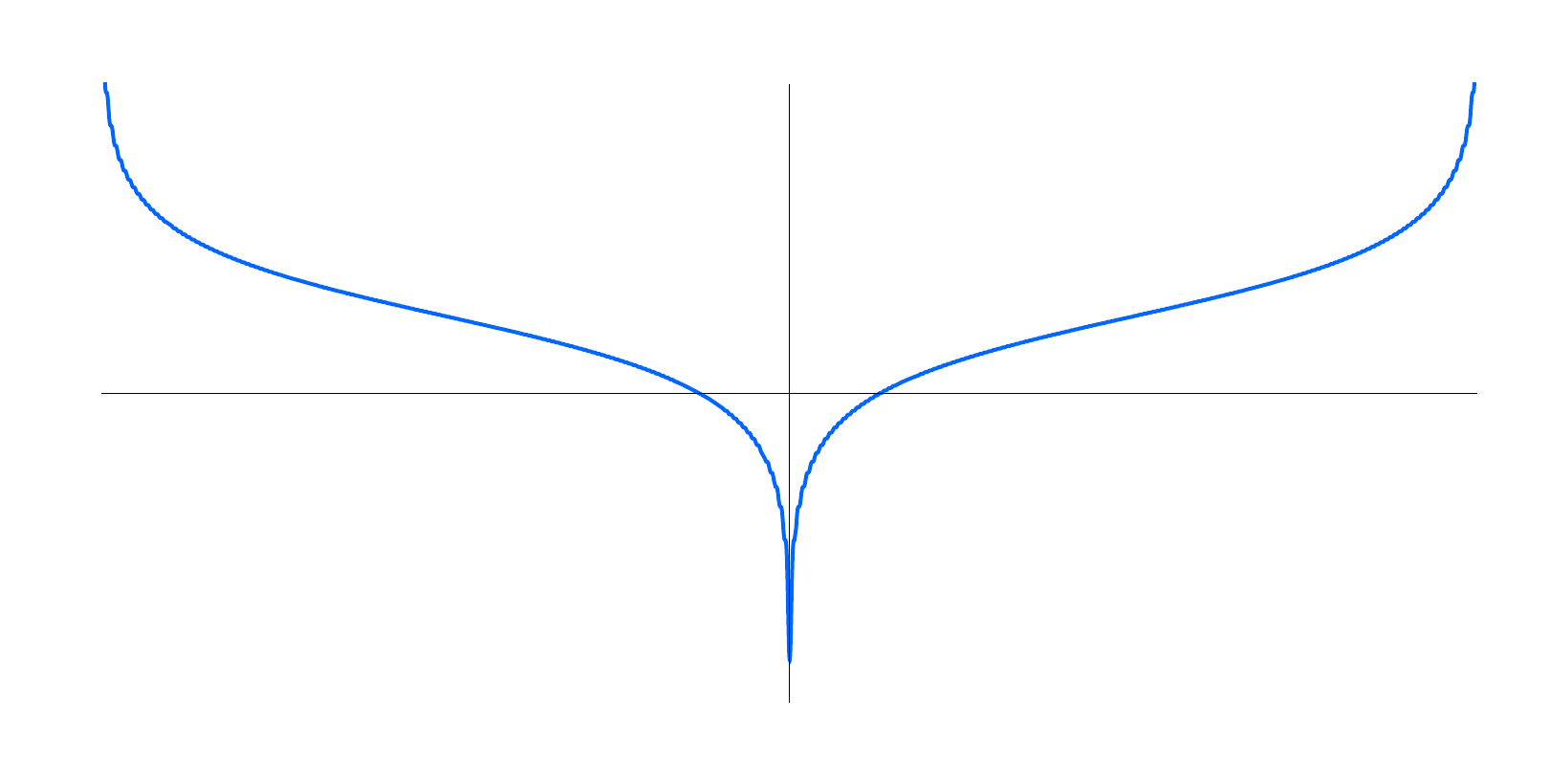}
    \subcaption*{$t=\frac \pi 2$}
           \end{subfigure}
              \begin{subfigure}{.5\textwidth}
\centering
       \includegraphics[height=1.2in]{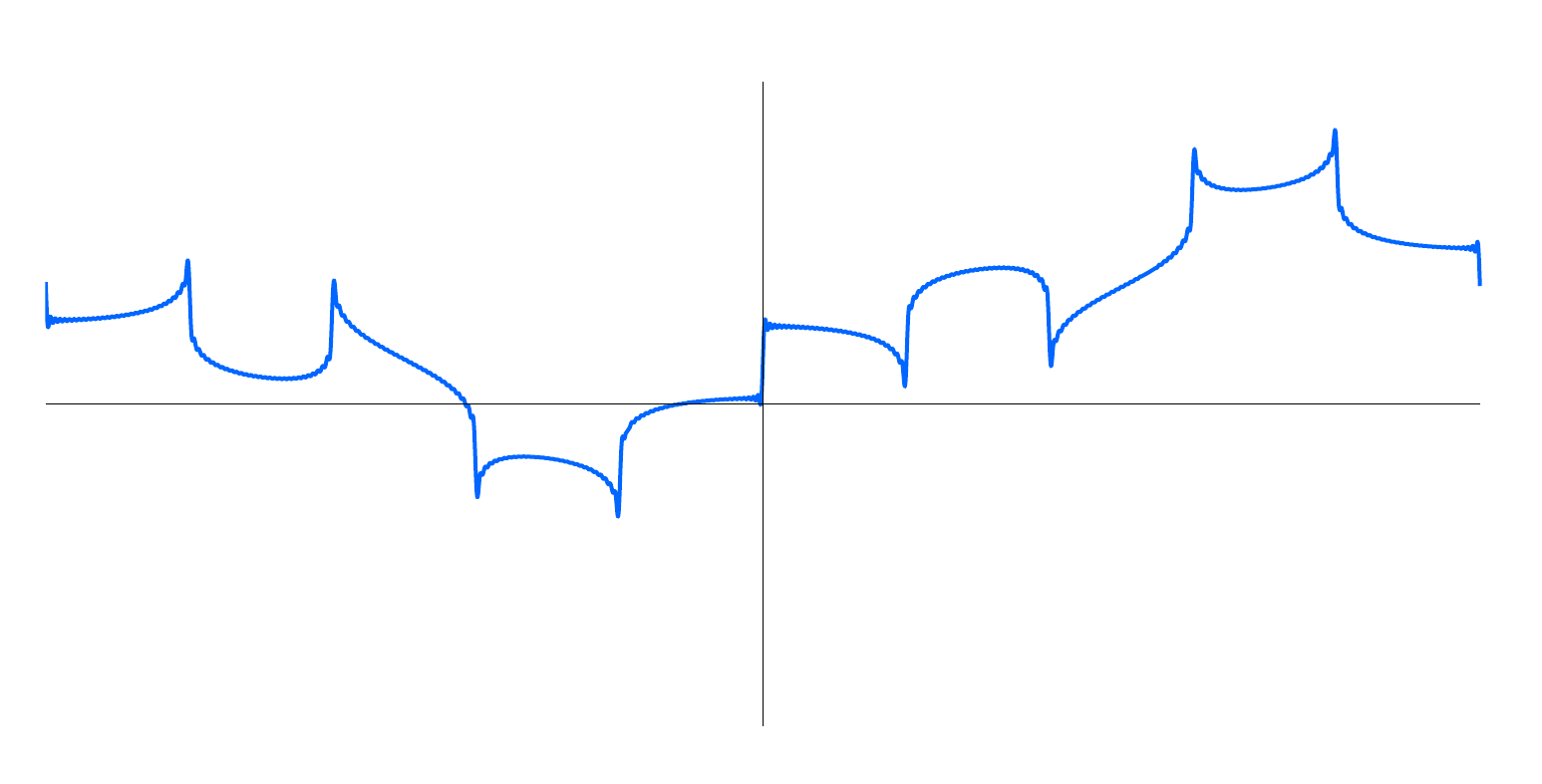}
    \subcaption*{$t=\frac \pi 5$}
           \end{subfigure}
              \begin{subfigure}{.5\textwidth}
\centering
       \includegraphics[height=1.2in]{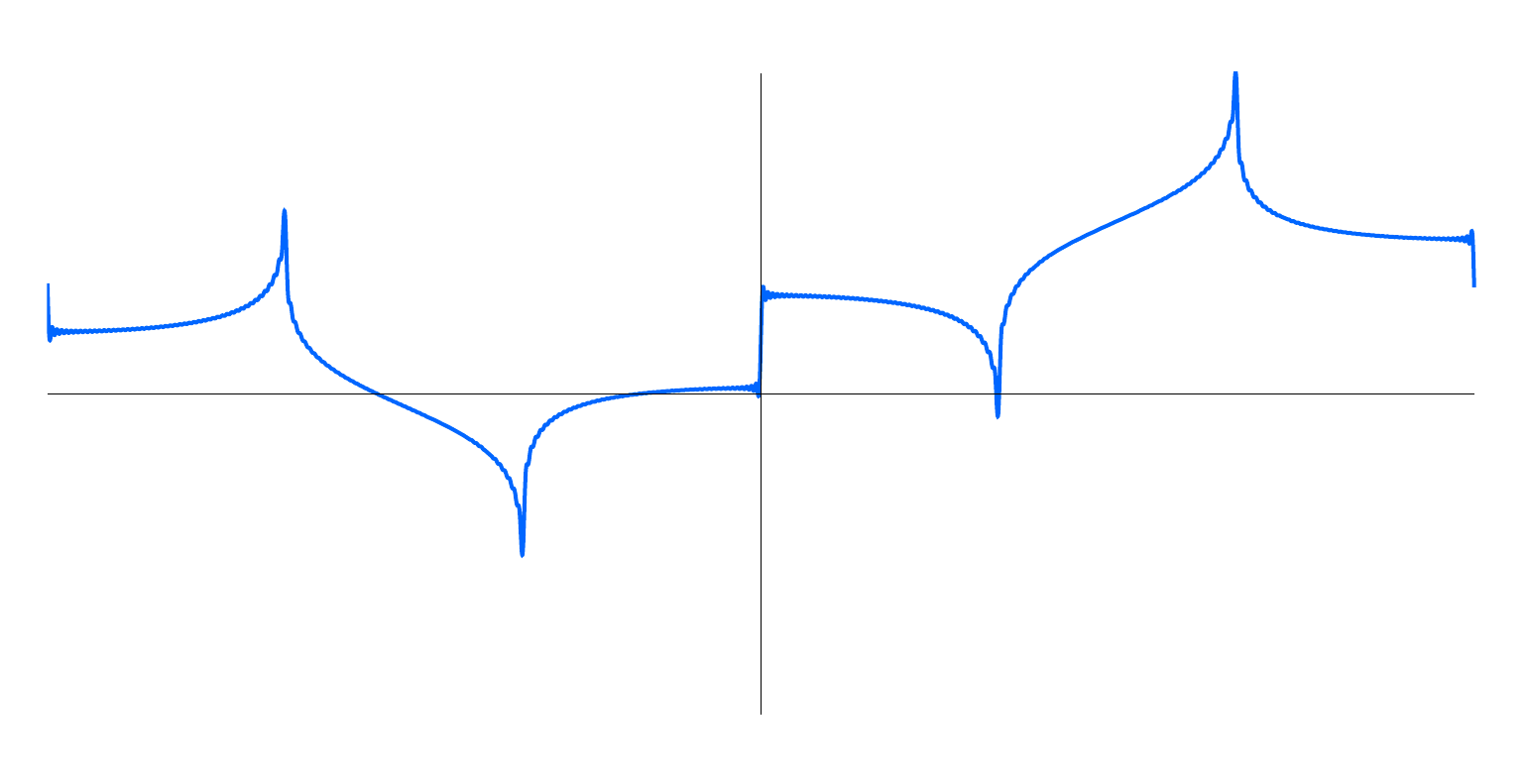}
    \subcaption*{$t=\frac \pi 6$}
           \end{subfigure}
           \hspace{.2in}
               \begin{subfigure}{.5\textwidth}
\centering
      \includegraphics[height=1.2in]{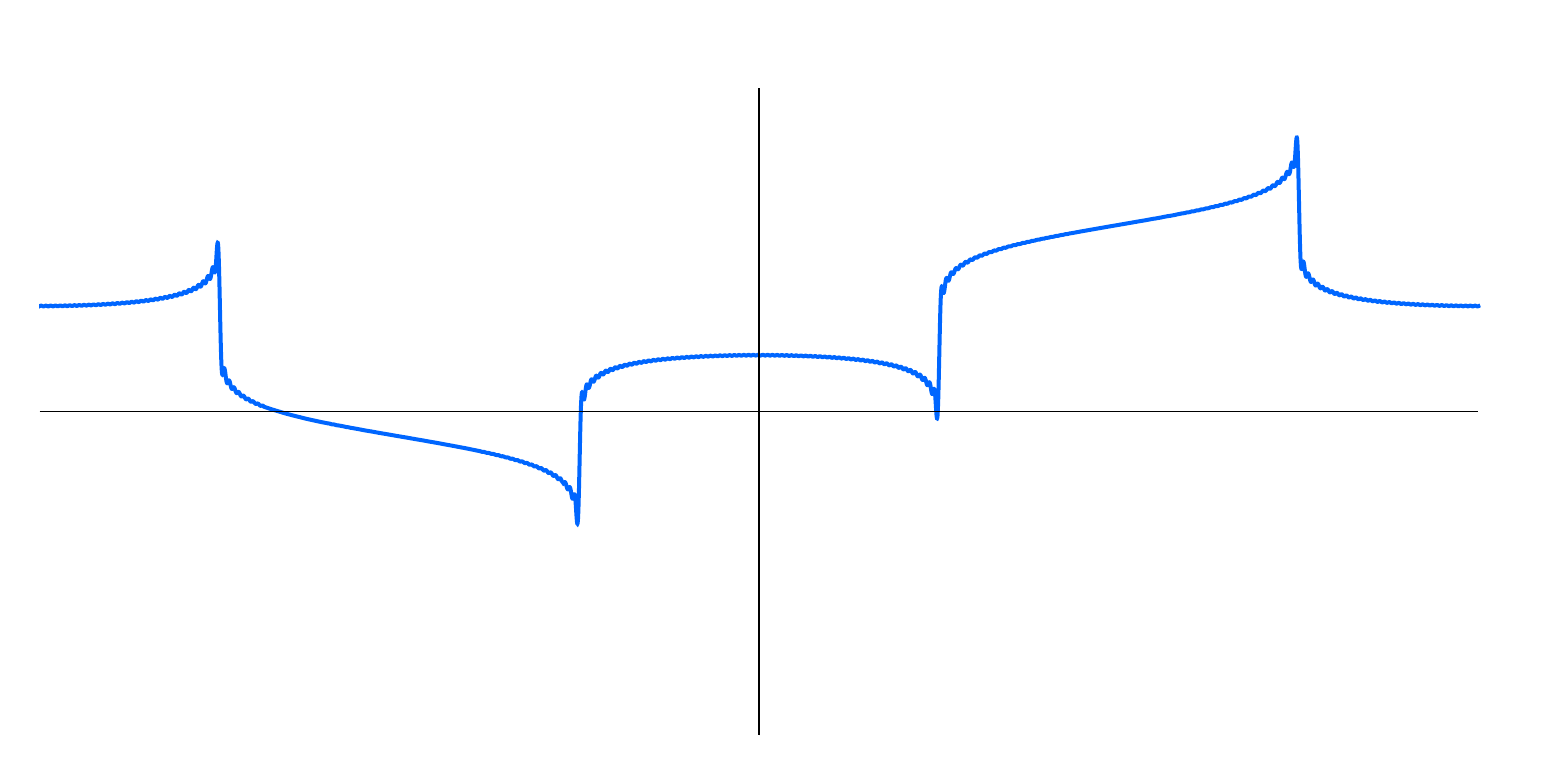}
    \subcaption*{$t=\frac \pi 8$}
           \end{subfigure}
               \hspace{.2in}
               \begin{subfigure}{.5\textwidth}
\centering
       \includegraphics[height=1.2in]{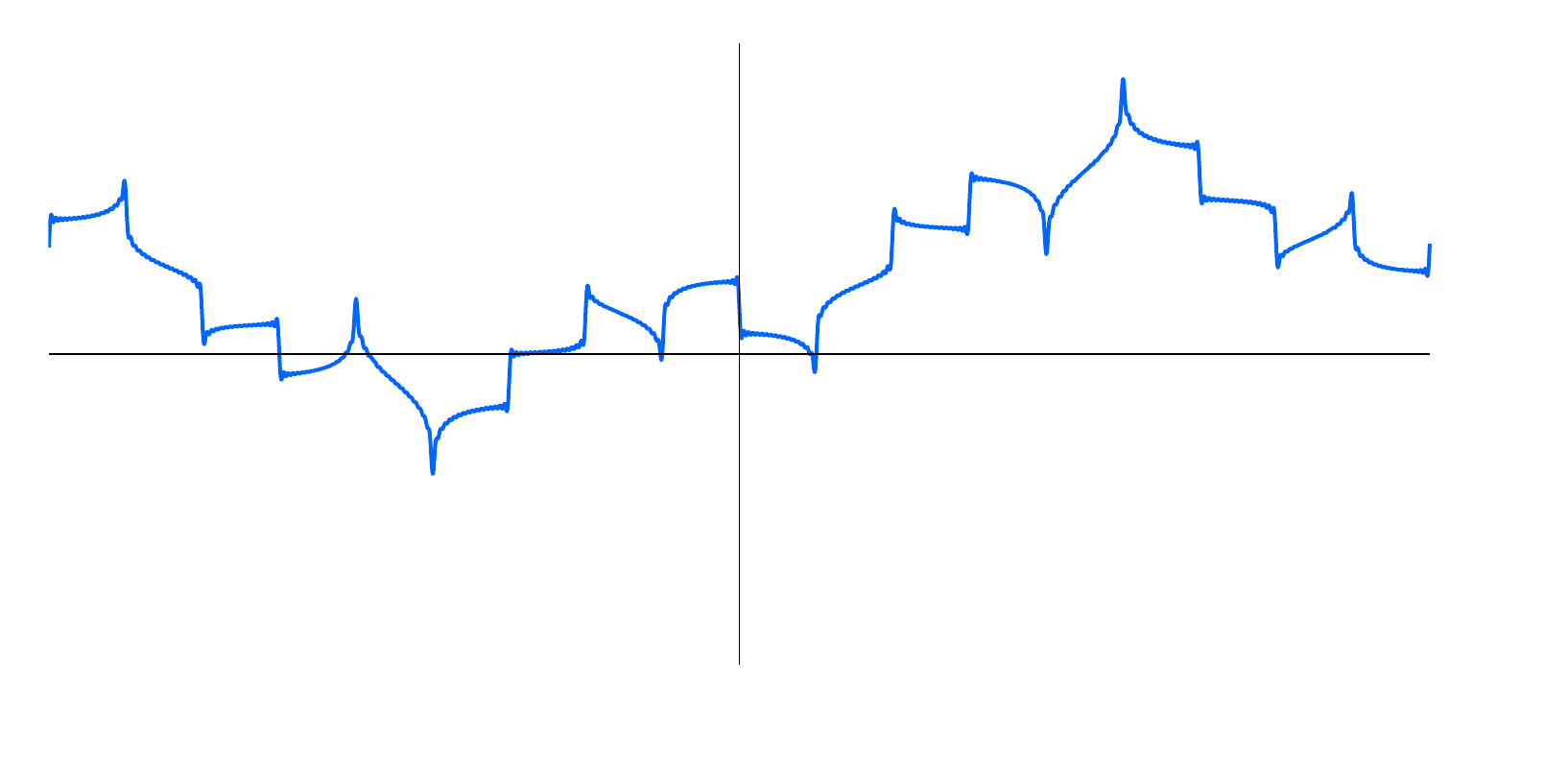}
    \subcaption*{$t=\frac \pi 9$}
           \end{subfigure}
        \caption{The solution of the Riemann problem for linearised BO equation at times that are rational multiples of $\pii$. The horizontal  axis is from $-\pii$ to $\pii$, the vertical axis is from $-2$ to $2$.  } 
        \label{fig:BOfig}
    \end{figure}

The solution of the Riemann problem  
at the rational time $t = p\pii/q$ can be written as
\Eq{lBOstep}
$$\eeq{u\Pa{\frac pq\,\pi,x} = \fra 2 + \frac 2\pi \Sumoi j\frac{\sin\bbk{(2\:j+1)(x - (2\:j+1)\:p\pii/q)}}{2\:j+1}
\\=\fra 2 + \frac 2\pi  \Sumo j{q-1} \Sumoi n\frac{\sin\bbk{(2\:n\: q + 2\:j+1)(x - (2\:n\: q + 2\:j+1)\:p\pii/q)}}{2\:n\: q + 2\:j+1}
\\= \fra 2 + \frac 2\pi \Sumo j{q-1} \Sumoi n\frac{\sin\bbk{(2\:n\: q + 2\:j+1)(x - (2\:j+1)\:p\pii/q)}}{2\:n\:q + 2\:j+1}
\\= \fra 2 + \frac 2\pi \Sumo j{q-1} S^{2q}_{2j+1}\Paz{x - \frac {(2\:j+1)\:p}q\pii},}$$
where we used \eq{SC}, \eq{E}, for the final equality.
Thus, apart from the constant term, the solution profile is a linear superposition of translated trigonometric hypergeometric functions \eq{Sjk1}.  
See Figure \ref{fig:BOfig} for some typical solution profiles. Hence, the cusps, which were first observed in \cite{COdisp} but not understood, now have a rigorous explanation.  In particular, we have shown that the solution profiles observed in \rf{COdisp} are logarithmic cusps of infinite height, a fact that was not suspected by the authors at the time, due to the graphics issues noted above.

We wish to obtain the fundamental solution $u = F(t,x)$ for the periodic linearised Benjamin--Ono equation, \ie the solution having the Dirac delta function as initial condition. To this end, we cannot simply manipulate the Fourier series of the Dirac delta due to its weak convergence, and we cannot just differentiate the step function solution since this produces a pair of Dirac deltas, and it is not so clear how to split off a single component.  Instead, we begin with the initial condition
\Eq{intdelta}
$$\req{u(0,x) = \gamma (x), \quad \mbox{where}\quad \gamma(x) = -\,\frac x{2\pii} + \frac{\sign x}2 = \fra \pi \, S^1_1(x)= \fra \pi \, \Sl_1(x), \\-\pii < x < \pi.}$$
We call $\gamma(x)$ the \is{integrated delta function}, since 
\Eq{dintdelta}
$$\gamma  '(x) = \dide (x) - \frac 1{2\pii}.$$
Its Fourier series is given by
\Eq{intdeltafs}
$$\gamma  (x) =
\fra\pi\Sumi k \frac{\sin k\:x}k.$$
The corresponding solution to the linearised Benjamin--Ono equation is
$$u(t,x) = \fra\pi\Sumi k \frac{\sin \bbk{k\:(x - k\: t)}}k.$$
At the rational time $t = p\pii/q$, this solution can be written as
$$\eeq{u\Paz{\frac pq\,\pi,x} =
\fra\pi\Sumi k \frac{\sin \bbk{k\:(x - k \:p\pii/q)}}k
= \Sum k{2\:q} \Sumoi n\frac{\sin\bbk{(2\:n\: q + k)(x - (2\:n\: q + k)\:p\pii/q)}}{2\:n\: q + k}
\\= \Sum k{2\:q} \Sumoi n\frac{\sin\bbk{(2\:n\: q + k)(x - k\:p\pii/q)}}{2\:n\:q + 2\:j+1}
= \Sum k{2\:q} S^{2q}_k\Paz{x - \frac {k\:p}q\pii}.}$$
Thus, in view of \eq{dintdelta}, the fundamental solution, i.e. the solution evolved from $\dide(x)=\gamma  '(x) + 1/({2\pii})$,  at rational times is formally obtained by differentiation:
\Eq{BOfundsol}
$$F\Paz{\frac pq\,\pi,x} = \frac 1{2\pii} + \pdd ux\Paz{\frac pq\,\pi,x} = \frac 1{2\pii} + \Sum k{2\:q} \od {S^{2q}_k}x \Paz{x - \frac {k\:p}q\pii},$$
with the explicit formula for the derivatives given in \pr{propSC}; see also equation \eq{dSkjl}.  The constant term in the fundamental solution \eqf{BOfundsol} is due  to the fact that the partial differential equation admits constant solutions.
It follows, recalling the expression \eq{dSkjl}, that at rational times the fundamental solution \eqf{BOfundsol} is a linear combination of delta functions combined with a linear combination of cotangents.  When convolved with a general initial condition, the delta components produce a standard revival phenomenon, being a linear combination of translates of the initial condition by integer multiples of $\pi/q$.  The cotangent terms are linear combinations of translates, by the same amount, of the  periodic Hilbert transform (\ref{perhtr}) of the initial condition. We have thus proved the following striking revival result.

\Th{lbof} At a rational time $\pi\:p/q$, the solution to the periodic initial-boundary value problem for the linearised Benjamin--Ono equation \eq{lBO} on the interval $-\pii < x < \pi$, with initial datum $u(0,x) = f(x)$ is a linear combination of translates $f(x + \pi\:j/q)$, for $\rgo j{2\:q-1}$,  of the initial condition and translates $g(x + \pi\:j/q)$ of its periodic Hilbert transform: $g(x) = \CH[f](x)$.

In particular, when the initial data is a step function, the first of the two components of~\eqref{Skjl} in the translated trigonometric polylogarithm solution \eqf{lBOstep} is comprised of piecewise constant translates of the step function\fnote{The individual summands are piecewise linear, but they combine to form a piecewise constant profile.}, while the logarithmic sine terms are translates of its circular Hilbert transform \eqref{perhtr}.


\begin{rem}
    A similar analysis applies to the case of the linearised KdV equation (\ref{lkdv}).  In that case, at the rational time $t = p\pii/q$, we can write the solution of the Riemann problem as
    $$\eeq{u\Pa{\frac pq\,\pi,x} = \fra 2 + \frac 2\pi \Sumoi j\frac{\sin\bbk{(2\:j+1)(x - (2\:j+1)^2\:p\pii/q)}}{2\:j+1}
    \\= \fra 2 + \frac 2\pi \Sumo j{q-1} \Sumoi n\frac{\sin\bbk{(2\:n\: q + 2\:j+1)(x - (2\:n\: q + 2\:j+1)^2\:p\pii/q)}}{2\:n\: q + 2\:j+1}
    \\= \fra 2 + \frac 2\pi \Sumo j{q-1} \Sumoi i\frac{\sin\bbk{(2\:n\: q + 2\:j+1)(x - (2\:j+1)^2\:p\pii/q)}}{2\:n\:q + 2\:j+1}
    \\= \fra 2 + \frac 2\pi \Sumo j{q-1} S^{2q}_{2j+1}\Paz{x - \frac {(2\:j+1)^2\:p}q\pii}.}$$
    However, it is known, \rf{Odq}, that the revivals in this case are all piecewise constant at rational times, without any cusps, as shown in Figure \ref{kdvfig}.  This is a consequence of the fact that in the sum \eq{Skjsum}  all the cusped Hilbert transform components cancel out.
\end{rem}

\subsection{Revival for the Linearised Intermediate Long Wave Equation}
\label{rILW}

We now turn to the periodic Riemann problem for the linearised Intermediate Long Wave (ILW) equation \eq{lILW} with convolution kernel \eq{Czeta} and dispersion relation \eq{ILWdr}.

 The solution to the Riemann problem can thus be written as 
$$
u(t,x)= \frac 1 2 +\frac 2 \pi \sum_{k=0}^\infty\frac{\sin\bbk{(2k+1)\bpa{t/\delta +x-(2k+1)t\coth[\delta(2k+1)]}}}
{2k+1}.
$$
Numerical experiments indicate that the solution of the Riemann problem for the  ILW equation, at least for  $\delta$ not too small,  presents a ``cusped'' behaviour at rational times analogous to that observed for the solution of the BO equation.  To investigate this phenomenon and its domain of validity, we start by comparing ILW with BO.  


It is well known that in the limit as $\delta \to \infty$, the solutions of ILW tend to the solutions of BO. This can be justified by considering  the asymptotic series 
\Eq{coth}
$$
\coth z \>\sim\> 1 + 2 \sum_{n=1}^\infty e^{-2nz} \>= \>1+O(\re^{-2z}), \qquad z \longrightarrow \infty .
$$
Thus, for $\delta>0$  and  $k\to \infty$, 
$$
\omega_{\delta}(k)=k^2\coth(\delta k)-\frac k \delta\ \sim \ k^2-\frac k \delta=k^2\Paz{1-\frac 1 {\delta k}}.
$$
To study the validity of this approximation, let
  $$\req{c_{BO} (k) = \frac{\omega_{BO} (k)}k = \abs k, \\ c_\delta (k) = \frac{\omega _\delta(k)}k = k \coth(\delta \, k) - \fra \delta ,}$$
be the corresponding phase velocities.
The solution to the periodic Riemann problem has the form
$$u_\delta(t,x) = \Sumii k b_k \exp\bbk{\i k\:\paz{x - c_\delta(k)\:t}},$$
where $b_k$ are the Fourier coefficients of the initial step data.
Note that for $|k| \gg 0$,
$$\Abs{k \coth(\delta \, k) - \abs k\ostrut84} = \abs k \frac{2 e^{-2\: \delta \: \abs k}}{1 -e^{-2\: \delta \: \abs k}} \leq a \,\abs k \,e^{-2\: \delta \: \abs k},$$
where the constant $a$ is close to $2$.  Thus, 
$$c_\delta (k)  =  c_{BO} (k) - \fra \delta + \varepsilon _\delta (k)$$
where
\Eq{eps}
$$\abs {\varepsilon _\delta (k)} \leq a \,\abs k \,e^{-2\: \delta \: \abs k}.$$
Now write
\be
u_\delta(t,x) = u_{BO}(t,x+ t/\delta ) + v(t,x).
\label{vdefwithu}\ee
The first term is just a phase shifted version of the BO solution, and hence cusped at rational times.
The difference can be expressed as a Fourier series:
\Eq{vfs}
$$v(t,x) = - \Sumii k b_k \exp\bbk{\i k( \abs k - 1/ \delta )\:t}\paz{1 - \exp\bbk{\i k\,\varepsilon _\delta (k)}} e^{\i k\:x}.$$
Now, using \eq{eps}, for $|k| \gg 0$,
$$\Abs{1 - \exp\bbk{\i k\,\varepsilon _\delta (k)}\ostrut84} \leq \atilde \,\abs k\,\varepsilon _\delta (k) \leq a \,\atilde \,k^2 \,e^{-2\: \delta \: \abs k},$$
for some $\atilde \approx 1$.
Thus, because $b_k$ are bounded (in fact decay like $1/k$), and, for fixed $t$, the second factor in the summand \eq{vfs} has modulus $1$, the Fourier series for $v(t,x)$ has coefficients that decay exponentially fast to $0$ as $\abs k \to \infty$.  This implies that  $v(t,x)$ is a $\Cci$ function.  Moreover, except for maybe its first few modes, one can estimate its $L^ \infty $ norm to be very small.  This implies that the  ILW solution differs from a phase shifted version of the BO solution by a smooth function which is 
close to a rational function of the first few modes. 
In turns, the ILW fundamental solution at $t=\pii p/q$ can be approximated by shifting the BO fundamental solution by $t/\delta$. See (\ref{ILWApp}).

         \begin{figure}
            \begin{subfigure}{1\textwidth}
       \centering
       \includegraphics[height=2in]{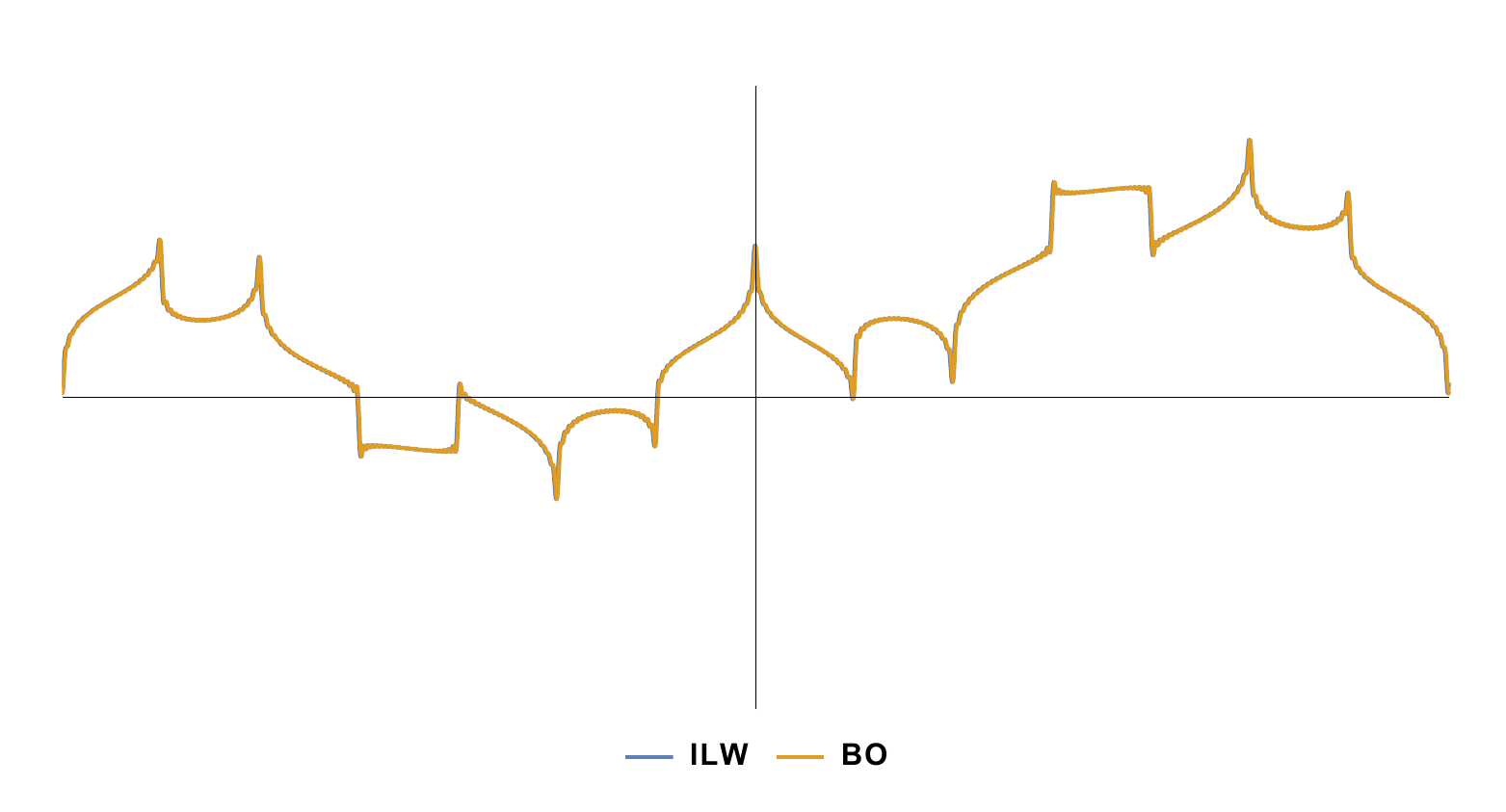}
    \subcaption*{ILW with $\delta=100$ vs BO }
           \end{subfigure}
             \caption{Solutions of the periodic Riemann problem for the linearised ILW and BO  at time $t=\pi/7$ compared for   $\delta=100$. The horizontal  axis is from $-\pii$ to $\pii$, the vertical axis is from $-2$ to $2$.  
             }
            \label{fig:ILWfig1}
           \end{figure}

\begin{figure}[h]
  \begin{subfigure}{1\textwidth}
        \centering
       \includegraphics[height=2in]{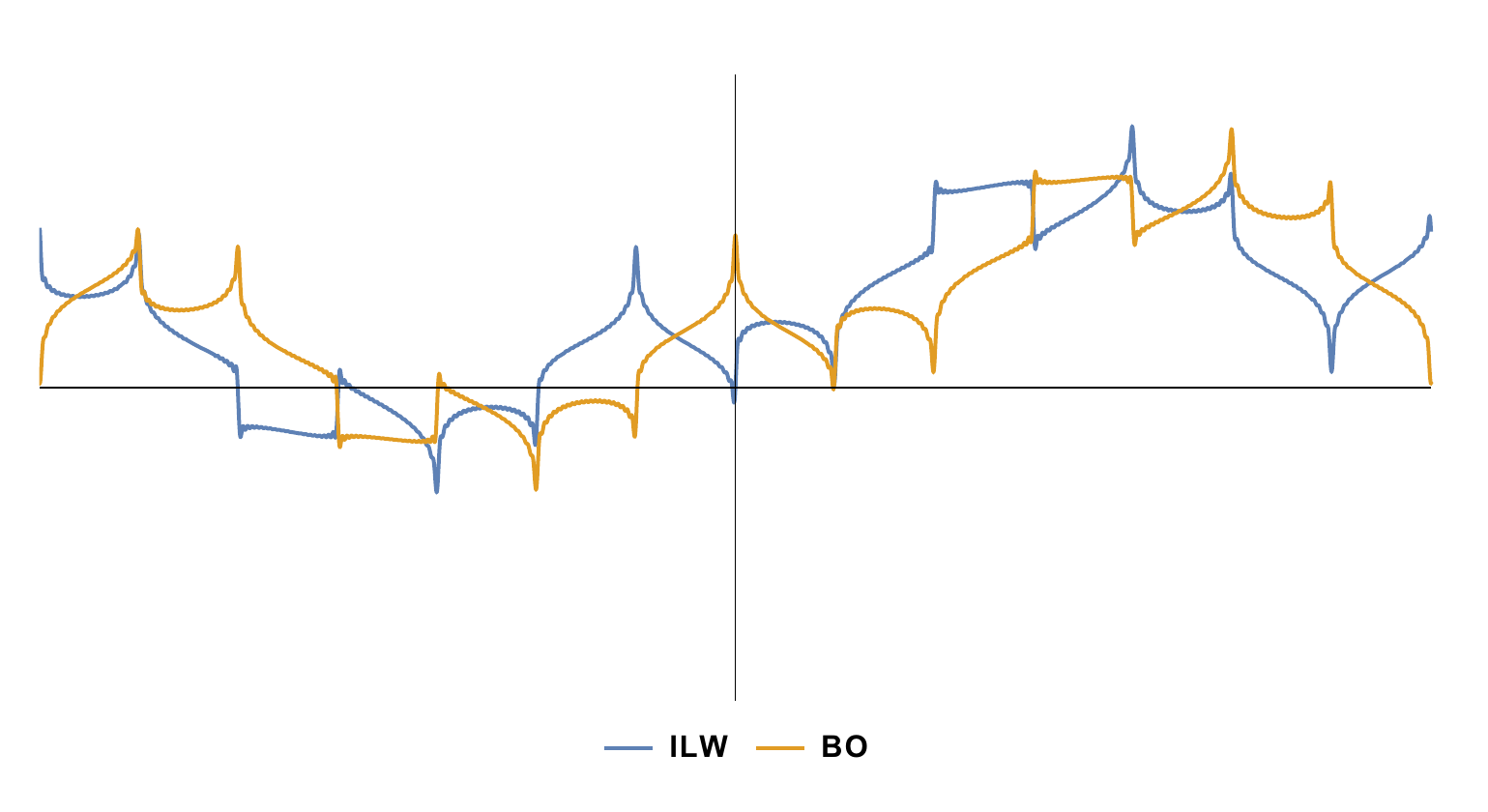}
       \subcaption*{ILW with $\delta=1$  vs BO, evaluated at the same value of $x$  - the shift by $t/\delta$ is in clear evidence}
           \end{subfigure}
             \begin{subfigure}{1\textwidth}
        \centering
       \includegraphics[height=2in]{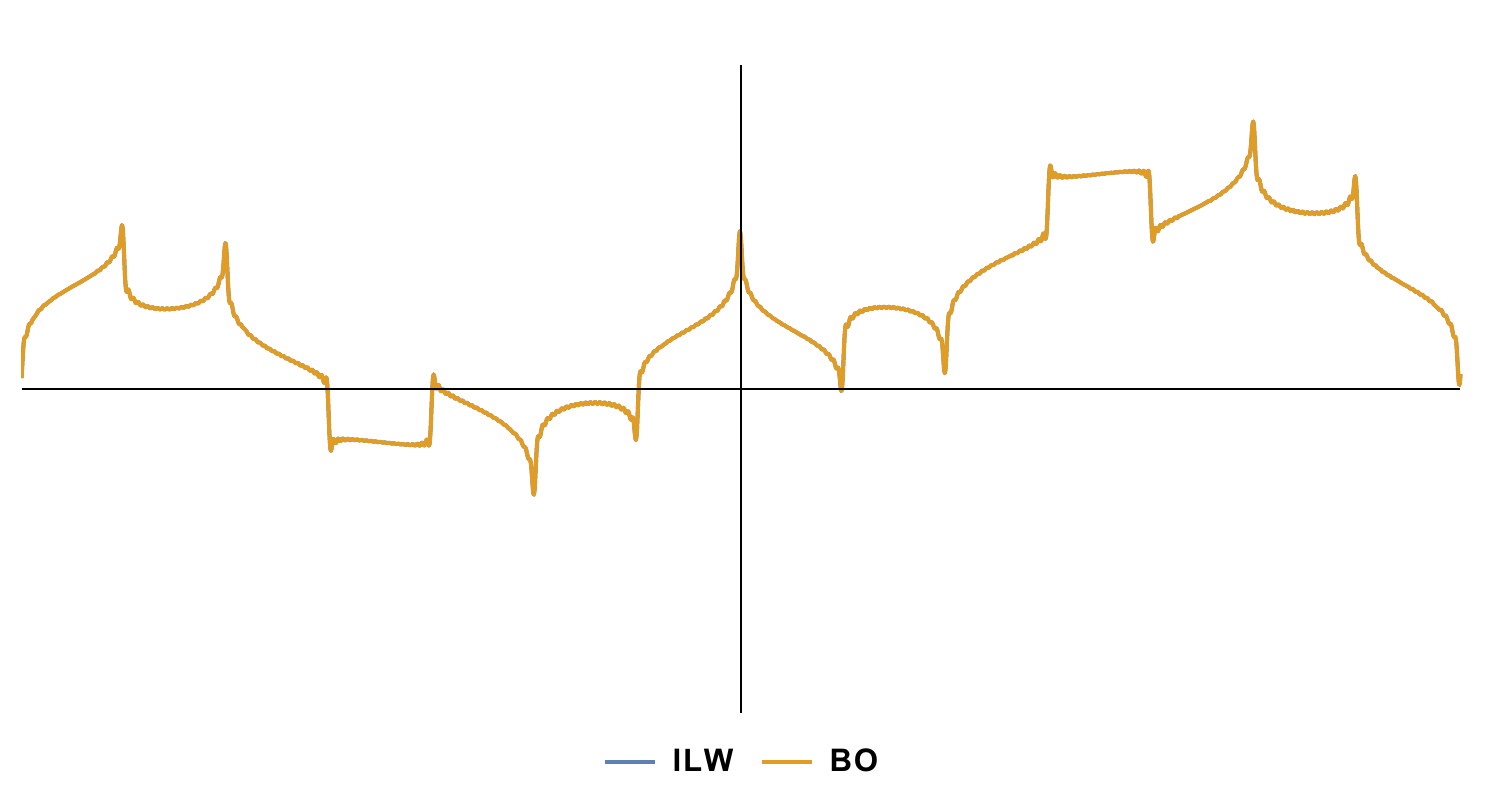}
       \subcaption*{ILW, with $\delta=1$, evaluated  at $x$ vs BO evaluated at $(t,x+t/\delta)$  }
           \end{subfigure}
            \caption{Solutions of the periodic Riemann problem for the linearised ILW with $\delta=1$ and BO  at time $t=\pi/7$. The solutions are evaluated at the same value of $x$ in the top graph, while BO includes a spatial shifted of $t/\delta$ in the bottom graph. The horizontal  axis is from $-\pii$ to $\pii$, the vertical axis is from $-2$ to $2$. }\label{fig:ILWfig2}
           \end{figure}

\smallskip
Numerically, we can test this observation for the associated Riemann problems. In Figures~\ref{fig:ILWfig1} and \ref{fig:ILWfig2}, we consider approximations obtained by truncating the series representation of $u_\delta$ and $u_\infty$. The experiment indicates that the phase shift depends on the ratio $t/\delta$, hence for fixed $\delta$ it increases with $t$, while the difference $v(t,x)$ depends on the size of $\delta$. Note that $\delta$ does not have to be chosen too large in order for the phenomenon to be observed, see Figure~\ref{fig:ILWfig2}. 

Using 
(\ref{vdefwithu}),
we can express
$u_\delta-v$ as 
\begin{align*}
(u_\delta-v)(t,x)&= \frac 1 2 +\frac 2 \pi \sum_{k=0}^\infty\frac{\sin\{(2k+1)[t/\delta+x-(2k+1)t)]\}}
{2k+1}\\&=\frac 1 2 +\frac 2 \pi \sum_{k=0}^\infty\frac{\sin\{(2k+1)[x-(2k+1)t+t/\delta)]\}}{2k+1}.
\end{align*}
Hence, at time $t=\pi p/q$ this gives 
 $$
(u_\delta-v)\Paz{\frac{p\pii}{q},x}= \frac 1 2 +\frac 2 \pi \sum_{k=0}^\infty\frac{\dsty\sin\Bk{(2k+1)\Pa{x-(2k+1)\frac{p\pii}{q}-\frac{p\pii}{q\delta}}}}
{2k+1}.
$$
Writing $k=nq+j$, $j=0,\ldots,q-1$, the terms $2nq\pi/q$ do not contribute to the evaluation of the sine function, so we can rewrite this sum as
$$
\ibeq{100}{\sum_{k=0}^\infty\frac{\dsty\sin\Bk{(2k+1)\Pa{x-(2k+1)\frac{p\pii}{q}-\frac{p\pii}{q\delta}}}}
{2k+1}\\=\sum_{j=0}^{q-1}\sum_{n=0}^\infty\frac{\dsty\sin\Bk{(2nq+2j+1)\Pa{x-(2j+1)\frac{p\pii}{q}+\frac{p\pii}{q\delta}}}}{2nq+2j+1},}
$$
hence
\be
(u_\delta-v)\Paz{\frac{p\pii}{q},x}=  \frac 1 2 +\frac 2 \pi\sum_{j=0}^{q-1}S_{2j+1}^{2q}\Bk{x-(2j+1)\frac{p\pii}{q}+ \frac{p\pii}{q\delta}},
\label{ILWApp}\ee
which is the explicit expression for (\ref{vdefwithu}).

\begin{figure}[t]
  \begin{subfigure}{1\textwidth}
       \centering
       \includegraphics[height=2in]{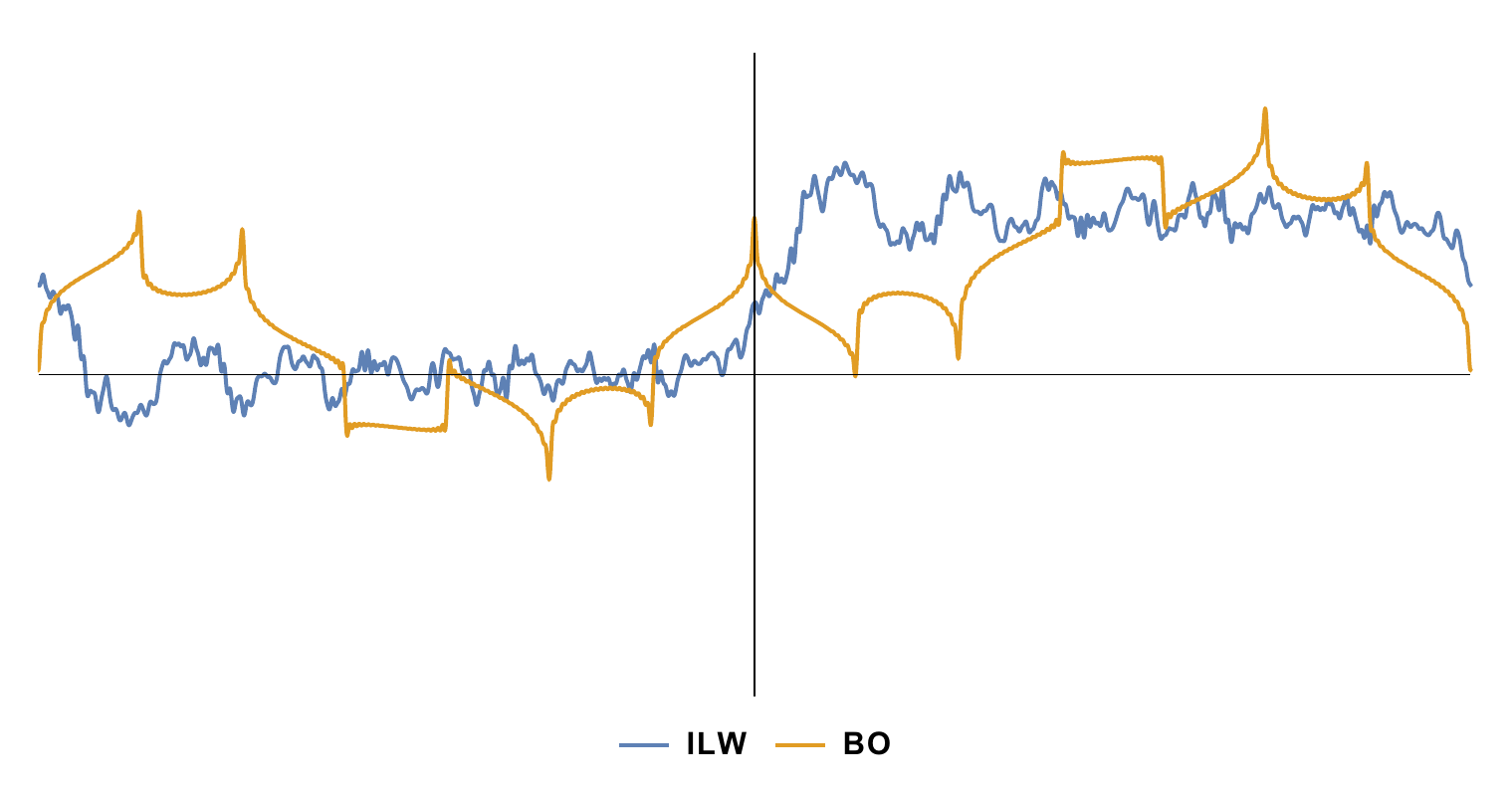}
       \subcaption*{ILW with $\delta=0.01$ vs BO}
           \end{subfigure}
            \begin{subfigure}{1\textwidth}
        \centering
       \includegraphics[height=2in]{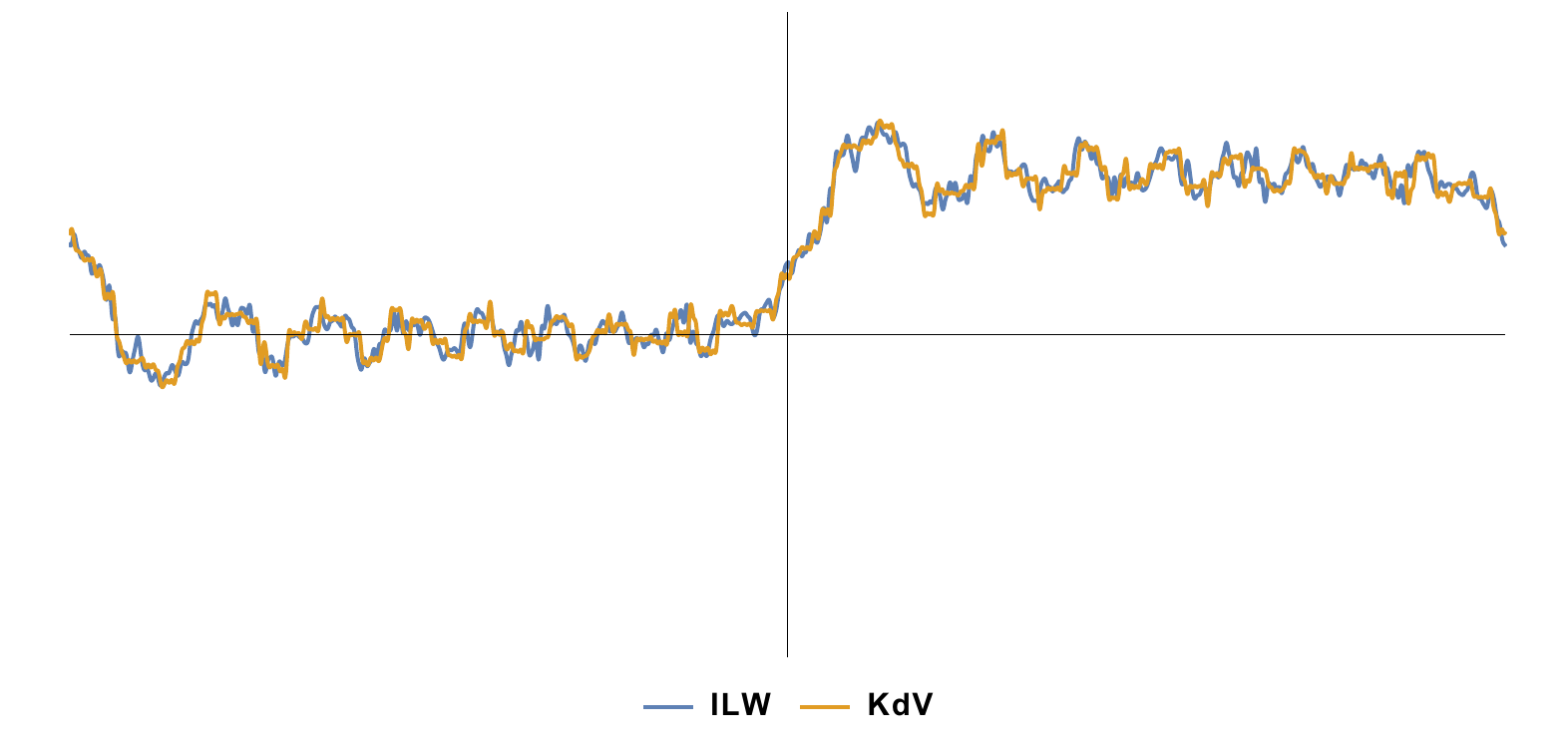}
    \subcaption*{ILW with $\delta=0.01$ vs KdV  }
           \end{subfigure}
             \caption{Solutions of the periodic Riemann problem for the linearised ILW, with  $\delta=.01$, compared with the solution of linearised BO and (rescaled) KdV at time $t=\pi/7$. The horizontal  axis is from $-\pii$ to $\pii$, the vertical axis is from $-2$ to $2$}
            \label{fig:ILWfig3}
           \end{figure}

In the final part of this section, let us investigate the asymptotic limit of the ILW equation as the depth parameter $\delta\to 0$.
The observed numerical behaviour illustrated above changes drastically for sufficiently small $\delta$, seemingly exhibiting revival phenomena more akin to those exhibited by the KdV equation. See Figure~\ref{fig:ILWfig3}.

As discussed in \cite{Saut}, if $u_\delta(t,x)$ denotes the solution of the ILW, then the function $v_\delta$ given by
\be
v_\delta(t,x)=\frac 3 \delta \,u_\delta\Paz{\frac 3 \delta t,x}, 
\label{kdvlim}\ee 
tends, as $\delta \to 0$ to the solution of the linearised KdV equation (\ref{lkdv}). 
Note that $v_\delta$ satisfies the rescaled equation
\be
(v_\delta)_t=\frac 3\delta {\cal L}[v_\delta] \roq{so that}
v_\delta(t,x)=\frac 1 {2\pi}\sum_{k=-\infty}^\infty\widehat v_{\delta}(0,k)\re^{\i kx- \delta\omega_\delta(k) t/3}.
\label{ilwres}\ee


Writing the dispersion relation as
$$	  \omega_\delta(k)=k^2\Bk{\coth(\delta k)-\frac 1 {\delta\, k}}$$
and using the Taylor expansion  around $\delta k=0$,
we can write
$$
\coth(\delta k)-\frac 1 {\delta k} =\frac 1 3\, \delta k-\frac 1 {45}\,(\delta k)^3+....,$$
so that
$$
\omega_\delta(k)\sim \f3\, \delta \,k^3+O(\delta^3k^5)= \f3\, \delta \, \omega_{kdv}(k)+O(\delta^3k^5),
$$
which implies that the expression (\ref{ilwres}) differs from the solution  of the linearised KdV equation by a multiplicative term of order 
$
\re^{-t\delta^3 k^5}$. For any finite value $k$, this term  tends to $1$  as $\delta\to 0$.

This calculation justifies the rescaling  (\ref{kdvlim}). 
Numerical calculations obtained summing a finite number of terms and choosing $\delta$ small  confirm the validity of the approximation above: the solution of the ILW at time $t$  and the solution of linearised KdV  at time $\delta t/3$, for small enough $\delta$,  are close - see figure \ref{fig:ILWfig3}.

\subsection{Revival for the Linearised Smith Equation}
\label{smith}

We now turn to the periodic Riemann problem for the linearised Smith equation \eq{Smith} with convolution kernel \eq{perSmithker}.
Let us begin by rewriting the dispersion relation \eq{Smithdr} in the form
$$
    \omega_S(k) = k\sqrt{\frac 1\delta+k^2} = k^2 + \fra{2\delta} \,f\left(\frac 1{\delta k^{2}}\right),
$$
for
$$
    f(z) =  \frac{2}{z}\,\bpa{ \sqrt{1+z} - 1 } = \sum_{n=0}^\infty z^n \frac{(-1)^n(2n-1)!}{(n+1)!(n-1)!2^{2n-1}},
$$
where the series is uniformly convergent for $z\in[0,1]$.


At time $t = p\pii/q$, we are interested in terms of the form
$$
    \frac{1}{k}\sin\left[ kx - t\,\omega_S(k) \right] = \frac{1}{k} \sin\left[ k\left(x-\frac{p\pii k}{q}\right) \right] \mathcal{C}\left(\frac{1}{k^2};\delta\right) - \frac{1}{k} \cos\left[ k\left(x-\frac{p\pii k}{q}\right) \right] \mathcal{S}\left(\frac{1}{k^2};\delta\right),
$$
where, using the Taylor expansion at $z=0$,
\begin{align*}
    \mathcal{C}(z;\delta) &= \cos\left(\frac{p\pii}{2q\delta}f\left(\frac{z}{\delta}\right)\right) = \mathcal{C}(0;\delta) + \mathcal{C}'(0;\delta)z + \frac{\mathcal{C}''(\zeta_{\mathcal{C}};\delta)}{2}z^2, \\
    \mathcal{S}(z;\delta) &= \cos\left(\frac{p\pii}{2q\delta}f\left(\frac{z}{\delta}\right)\right) = \mathcal{S}(0;\delta) + \mathcal{S}'(0;\delta)z + \frac{\mathcal{S}''(\zeta_{\mathcal{S}};\delta)}{2}z^2,
\end{align*}
for suitable $\zeta_C,\zeta_S$ with $\lvert\zeta_C\rvert,\lvert\zeta_S\rvert\leq\lvert z \rvert\leq 1$.
Subtituting and matching terms above, yields
\begin{align*}
    \mathcal{C}(z;\delta) &= \cos\left(\frac{p\pii}{2q\delta}\right) + \frac{p\pii}{8q\delta^2}\sin\left(\frac{p\pii}{2q\delta}\right)z + \varepsilon_\mathcal{C} z^2, \\
    \mathcal{S}(z;\delta) &= \sin\left(\frac{p\pii}{2q\delta}\right) - \frac{p\pii}{8q\delta^2}\cos\left(\frac{p\pii}{2q\delta}\right)z + \varepsilon_\mathcal{S} z^2,
\end{align*}
where
\begin{equation*}
    \lvert\varepsilon_\mathcal{C}\rvert,\lvert\varepsilon_\mathcal{S}\rvert \leq \varepsilon:= \frac{p\pii}{8q\delta^3}\left(1+\frac{p\pii}{8q\delta}\right).
\end{equation*}
Note that there exists an alternative bound on the $\varepsilon_\mathcal{S}$ with quartic decay in $\delta^{-1}$, but that bound does not apply to $\varepsilon_\mathcal{C}$.

 \begin{figure}[h]
  \begin{subfigure}{1\textwidth}
        \centering
       \includegraphics[height=2in]{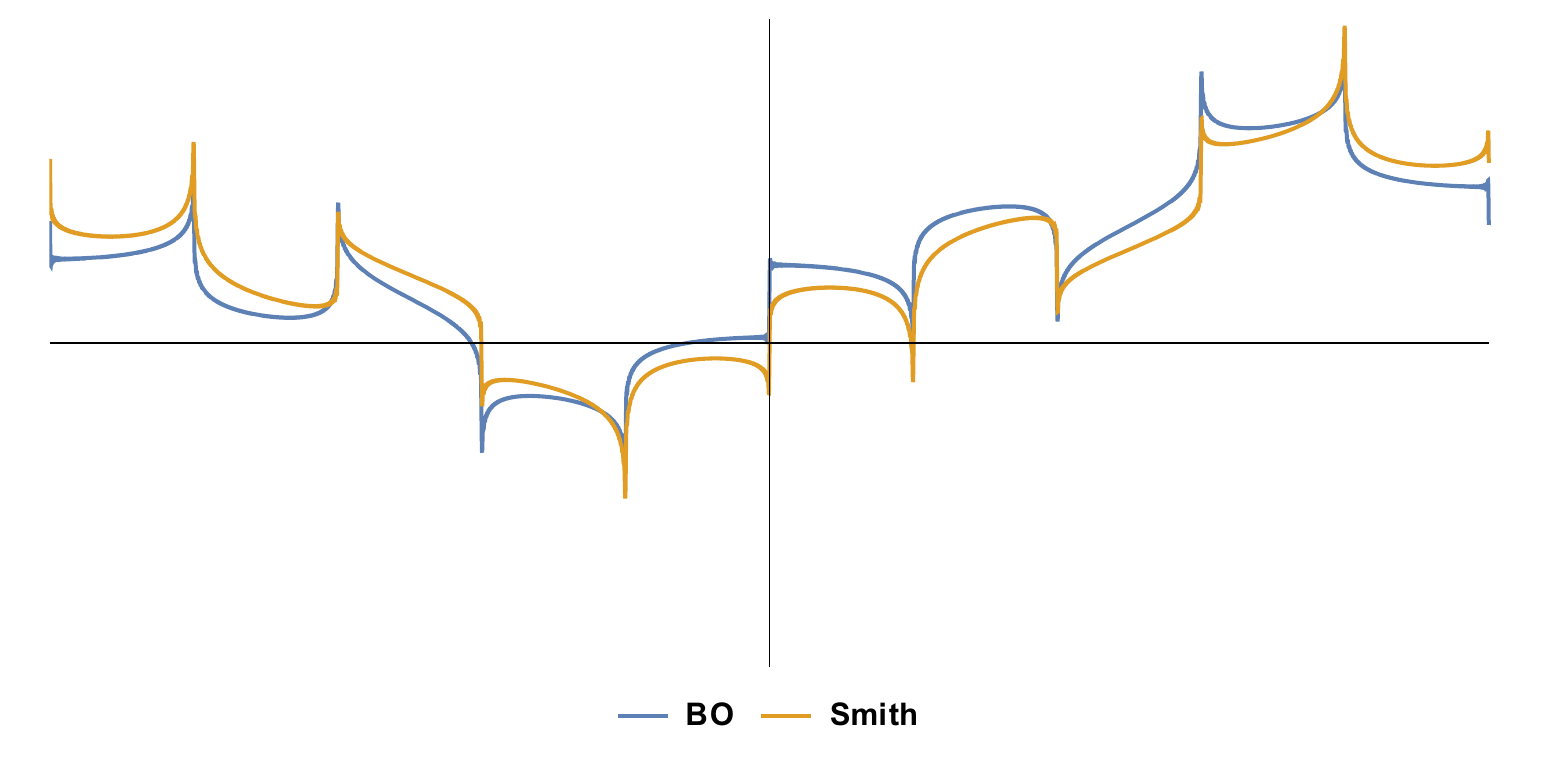}
       \subcaption*{Smith vs BO for $\delta=1$}
           \end{subfigure}
            \begin{subfigure}{1\textwidth}
        \centering
       \includegraphics[height=2in]{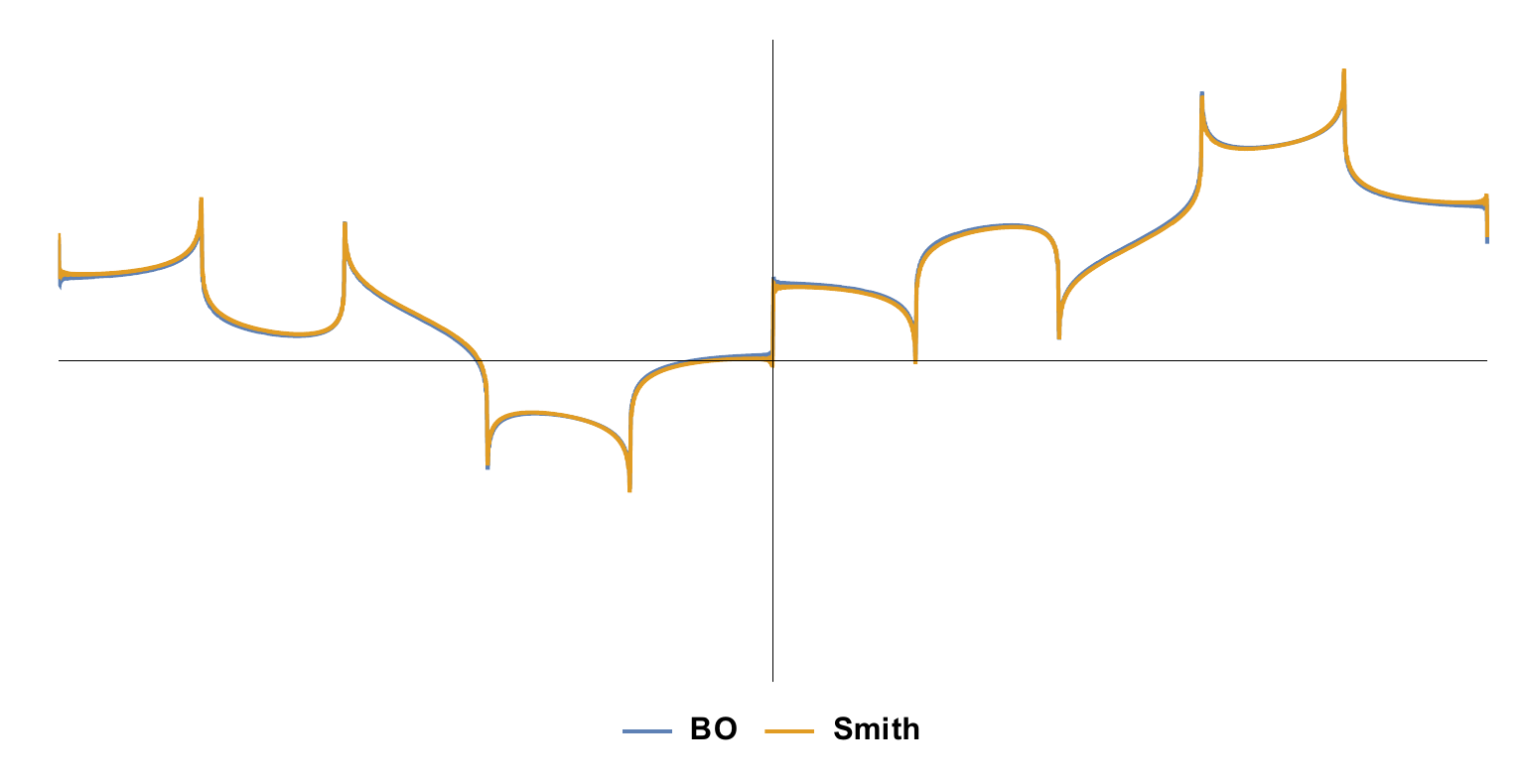}
    \subcaption*{Smith vs BO for $\delta=10$}
           \end{subfigure}
             \caption{Solutions of the periodic Riemann problem for the linearised Smith equation compared with the solution of linearised BO  at time $t=\pi/5$, for two values of $\delta$. The horizontal  axis is from $-\pii$ to $\pii$, the vertical axis is from $-2$ to $2$}
            \label{fig:ILWfig4}
           \end{figure}

Given that the dispersion relation for the Smith equation is an odd function, we can use the general \eqf{sol} and the explicit expression \eqref{eqn:SCdefn} to represent the solution $u(t,x)$ to the periodic Riemann problem for the Smith equation at time $t=p\pii/q$ as
\begingroup
\allowdisplaybreaks
    \begin{align*}
    u\left(\frac{p\pii}{q},x\right) &= \frac{1}{2} + \frac{2}{\pi} \sum_{j=0}^\infty \frac{\dsty\sin\left[(2j+1)\left(x-\frac{p(2j+1)\pi}{q}\right) - \frac{p\pii}{q\delta}f\left(\frac{1}{\delta(2j+1)^2}\right)\right]}{2j+1} \\
    &= \frac{1}{2} + \frac{2}{\pi} \sum_{j=0}^\infty \frac{\dsty\sin\left[(2j+1)\left(x-\frac{p(2j+1)\pi}{q}\right)\right]\mathcal{C}\left(\frac{1}{(2j+1)^2};\delta\right)}{2j+1} \\ &\hspace{8em}{}-{} \frac{2}{\pi} \sum_{j=0}^\infty \frac{\dsty\cos\left[(2j+1)\left(x-\frac{p(2j+1)\pi}{q}\right)\right]\mathcal{S}\left(\frac{1}{(2j+1)^2};\delta\right)}{2j+1}  \\
    &= \frac{1}{2} + \frac{2}{\pi}\cos\left(\frac{p\pii}{2q\delta}\right) \sum_{j=0}^{q-1} S^{2q}_{2j+1,1}\left(x-\frac{(2j+1)p\pii}{q}\right) \\ &\hspace{4em}{}-{} \frac{2}{\pi}\sin\left(\frac{p\pii}{2q\delta}\right) \sum_{j=0}^{q-1} C^{2q}_{2j+1,1}\left(x-\frac{(2j+1)p\pii}{q}\right) \\ &\hspace{4em}{}+{} \frac{p}{4q\delta^2}\cos\left(\frac{p\pii}{2q\delta}\right) \sum_{j=0}^{q-1} C^{2q}_{2j+1,3}\left(x-\frac{(2j+1)p\pii}{q}\right) \\ &\hspace{4em}{}+{} \frac{p}{4q\delta^2}\sin\left(\frac{p\pii}{2q\delta}\right) \sum_{j=0}^{q-1} S^{2q}_{2j+1,3}\left(x-\frac{(2j+1)p\pii}{q}\right) + \mbox{error}(x),
\end{align*}
\endgroup
in which
\begin{align*}
    \lvert\mbox{error}(x)\rvert
    &\leq \frac{2}{\pi} \left( \lvert\varepsilon_{\mathcal{C}} \rvert \left\lvert \sum_{j=0}^\infty \frac{\dsty\sin\left[(2j+1)\left(x-\frac{p(2j+1)\pi}{q}\right)\right]}{(2j+1)^5} \right\rvert \right. \\ &\hspace{10em} \left. + \lvert\varepsilon_{\mathcal{S}} \rvert \left\lvert \sum_{j=0}^\infty \frac{\dsty\cos\left[(2j+1)\left(x-\frac{p(2j+1)\pi}{q}\right)\right]}{(2j+1)^5} \right\rvert \right) \\
    &\leq \frac{4\varepsilon}{\pi} \sum_{j=0}^\infty \frac{1}{(2j+1)^5}
    < \frac{4\varepsilon}{\pi}\,\zeta(5) = \frac{p\,\zeta(5)}{2q\delta^3}\left(1+\frac{p\pii}{8q\delta}\right),
\end{align*}
where $\zeta(z)$ is the Riemann zeta function. 
Because the Fourier coefficients decay $\mathcal{O}(n^{-5})$, it follows that $\mbox{error}\in\ro C^{3}_{\mathrm{per}}[-\pii,\pi]$.
Therefore $u(p\pi/q,\cdot)$ is at least $\ro C^3$ away from the nodes.

Thus an approximate fundamental solution, at rational times, is obtained by differentiation: 
\begin{align*}
    F(t,x) &= \frac{2}{\pi} + \frac{\partial u}{\partial x}\left(\frac{p\pii}{q},x\right) \\
    &= \frac{2}{\pi} + \frac{1}{\pi}\cos\left(\frac{p\pii}{2q\delta}\right) \sum_{k=1}^{2q} \frac{dS_{k,1}^{2q}}{dx}\left(x-\frac{kp\pii}{q}\right)  {}-{} \frac{1}{\pi}\sin\left(\frac{p\pii}{2q\delta}\right) \sum_{k=1}^{2q} \frac{dC_{k,1}^{2q}}{dx}\left(x-\frac{kp\pii}{q}\right) \\ &{}-{} \frac{p}{8q\delta^2}\cos\left(\frac{p\pii}{2q\delta}\right) \sum_{k=1}^{2q} S_{k,2}^{2q}\left(x-\frac{kp\pii}{q}\right) {}+{} \frac{p}{8q\delta^2}\sin\left(\frac{p\pii}{2q\delta}\right) \sum_{k=1}^{2q} C_{k,2}^{2q}\left(x-\frac{kp\pii}{q}\right)
    \\&+\mbox{Error}(x),
\end{align*}
where the derivatives of the trigonometric polylogarithms are given in the appendix and the error term belongs to $\ro C^{2}_{\mathrm{per}}[-\pi,\pi]$.

We therefore conclude that, at rational times, the fundamental solution for the Smith equation  is a linear combination of  trigonometric polylogarithms up to a small error. Moreover, when the polylogarithms are convolved with the initial condition, a finite number of cusps arise. 


\section{Conclusions}

In this paper, we have identified and studied a novel form of the revival phenomenon. In the models analysed in the literature until now, the revival phenomenon for a spatially periodic PDE involves only translates of the initial data. In the new form we present here, the solution at rational times involves translates of both the initial data and its Hilbert transform. We have rigorously established and characterised this property for the linearised BO equation with piecewise constant initial condition. In this case, the Hilbert transform component,  at rational times, produces logarithmic cusps.   We have also shown that, both numerically and approximately,  the linearised ILW and Smith equations also exhibit similar phenomena. It is striking that, according to the numerics in \rf{COdisp}, at irrational times the graph of the solution appears to be a continuous but non-differentiable fractal curve, while at rational times it has infinite cusp singularities.  The alternation between these two scenarios mirrors the observed dispersive quantisation/fractalisation phenomena for equations with polynomial dispersion relations, like the linearised Schr\"odinger and KdV equations, for which the revival appears only in the form of a finite superposition of translates of the initial data.  Our analysis has relied on the properties of a novel class of special functions, the trigonometric polylogarithms and hypergeometric functions.

It would be interesting to investigate  further the dependence of this phenomenon on the periodicity of the problem. As is the case for equations with polynomial dispersion, we expect that the behaviour described in this paper persists, at least for boundary conditions that can be reformulated as periodic, e.g. Dirichlet, Neumann, or some pseudo-periodic conditions, as described in \cite{OSS} and \cite{Farmakis}. 

Another important research question arising from this work concerns the robustness of this novel revival phenomenon. We expect that it will persist in the nonlinear regime, specifically for the nonlinear BO, ILW, and Smith equations. In order to answer this question, or at least obtain numerical evidence, it will be necessary to devise suitably accurate numerical integration schemes. The numerical detection of the cusps will be particularly challenging.  These nonlinear equations  were derived in order to model waves in physical fluids.  However, unlike the optics and quantum mechanics experiments that confirm the revival and Talbot phenomena for the linear Schr\"odinger equation, \cite{BMS, VVS, YeaStr}, we expect that viscosity and other physical effects will make an experimental verification of these results very difficult. Hence numerical schemes able to accurately represent the revival phenomenon will be the only viable path to obtain evidence for the phenomenon in the nonlinear case.  



\section*{Acknowledgements}
The authors gratefully acknowledge support from Yale-NUS College workshop grant IG18-CW003, which funded all of them to attend a week-long workshop where the results presented here were first discovered and discussed.

\References

\key ABFS \paper Abdelouhab, L., Bona, J.L., Felland, M., Saut, J.C.; Nonlocal models for nonlinear, dispersive waves; Physica D: Nonlinear Phenomena; 40 (1989) 360-392\par
\key AFSS \paper Ablowitz, M.J., Fokas, A.S., Satsuma, J., Segur, H.; On the periodic intermediate long wave equation; J. Phys. A; 15 (1982) 781--786\par
\key Benjamin \paper Benjamin, T.B.; Internal waves of permanent form in fluids of great depth; J. Fluid Mech.; 29 (1967) 559--592\par
\key Berry \paper Berry, M.V.; Quantum fractals in boxes; J. Phys. A; 29 (1996) 6617--6629\par
\key BerryKlein \paper Berry, M.V., Klein, S.; Integer, fractional and fractal Talbot effects; J. Mod. Optics; 43 (1996) 2139--2164\par
\key BMS \paper Berry, M.V., Marzoli, I., Schleich, W.; Quantum carpets, carpets of light; Physics World; 14{\rm (6)} (2001) 39--44\par
\key BrezisPDEBook \book Brezis, H.; Functional Analysis, Sobolev Spaces and Partial Differential Equations; Springer, New York, 2010\par
\key BuNe \book Butzer, P.L., Nessel, R.J.; Fourier Analysis and Approximation (Vol I); Academic Press, New York and London, 1971\par
\key Calogero \paper Calogero, F.; Solution of the one-dimensional $n$-body problems with quadratic and/or inversely quadratic pair potentials; J. Math. Phys.; 12  (1971) 419--436\par
\key COdisp \paper Chen, G., Olver, P.J.; Dispersion of discontinuous periodic waves; Proc. Roy. Soc. London A; 469 (2013) 20120407--\par
\key ChErTz \paper Chousionis, V., Erdo\utxt gan, M.B., Tzirakis, N.; Fractal solutions of linear and nonlinear dispersive partial differential equations; Proc. London Math. Soc.; 110 (2015) 543--564\par
\key EMOT \book Erd{\'e}lyi, A., Magnus W., Oberhettinger, F. and Tricomi, F. G.; Higher Transcendental Functions. Vol. I; McGraw-Hill Book Company, New York-Toronto-London, 1953\par 
\key ErSh \paper Erdo\utxt gan, M.B., Shakan, G.; Fractal solutions of dispersive partial differential equations on the torus; Selecta Math.; 25 (2019) 11\par
\key ET \book Erdo\utxt gan, M.B., Tzirakis, N.; Dispersive Partial Differential Equations: Wellposedness and Applications; London Math. Soc. Student Texts, vol. 86, Cambridge University Press, Cambridge, 2016\par
\key Farmakis \paper Farmakis, G.; Revivals in Airy's   equation with quasi-periodic boundary conditions; preprint; (2020)\par
\key KS \paper Hofmanov{\'a}, M., Schratz, K.; An exponential-type integrator for the KdV equation; Numerische Mathematik; 136 (2017) 1117--1137\par
\key KS1 \paper Knöller, M.,  Ostermann, A., Schratz, K.; A Fourier Integrator for the Cubic Nonlinear Schrödinger Equation with Rough Initial Data; SIAM Journal on Numerical Analysis; 57 (2019) 1967--1986\par
\key Joseph \paper Joseph, R.I.; Solitary waves in a finite depth fluid; J. Phys. A; 10 (1977) L225--L228\par
\key King \book King, F.W.; Hilbert Transforms; Encyclopedia of Math. Appl. vols  124--125; Cambridge University Press, Cambridge, 2009\par
\key KKD \paper Kubota, T., Ko, D.R.S., Dobbs, L.D.; Weakly-nonlinear, long internal gravity waves in stratified fluids of finite depth; J. Hydronautics; 4 (1978) 157--165\par
\key Lew1981a \book Lewin, L.; Polylogarithms and Associated Functions; Elsevier North Holland, New York, 1981\par
\key MCK \paper McKean, H. P.; Boussinesq's equation on the circle; Commun. Pure Appl. Math.; 34 (1981) 599--691\par
\key MoserH \paper Moser, J.; Three integrable Hamiltonian systems connected with isospectral deformations; Adv. Math.; 16 (1975) 197--220\par
\key DLMF \book Olver, F.W.J., Lozier, D.W., Boisvert, R.F., Clark, C.W.{, eds.}; NIST Handbook of Mathematical Functions; Cambridge University Press, Cambridge, 2010\par
\key Odq \paper Olver, P.J.; Dispersive quantization; Amer. Math. Monthly; 117 (2010) 599--610\par
\key OSS \paper Olver, P.J., Sheils, N.E., Smith, D.A.; Revivals and fractalisation in the linear free space Schr\"odinger equation; Quart. Appl. Math.; 78 (2020) 161--192\par
\key OT \paper Olver, P.J., Tsatis, E.; Points of constancy of the periodic linearised Korteweg--deVries equation; Proc. Roy. Soc. London A; 474 (2018) 20180160--\par
\key Ono \paper Ono, H.; Algebraic solitary waves in stratified fluids; J. Phys. Soc. Japan; 39 (1975) 1082--1091\par
\key OskolkovV \inbook Oskolkov, K.I.; A class of I.M. Vinogradov's series and its applications in harmonic analysis; Progress in Approximation Theory; A.A. Gonchar and E.B. Saff, eds., Springer Ser. Comput. Math., 19, Springer, New York, 1992, pp. 353--402\par
\key DP \paper Pelloni, B, Dougalis, V.A.;Numerical solution of some nonlocal, nonlinear dispersive wave equations; Journal of Nonlinear Science; 10 (2000) 1-22\par
\key PBM3 \book Prudnikov, A.P., Brychkov, Y.A., Marichev, O.I.; Integrals and Series. Volume 3. More Special Functions; Gordon and Breach Sci. Publ., New York, 1990\par
\key Saut \preprint Saut, J.-C.; Benjamin-Ono and Intermediate Long Wave equations: Modeling, IST and PDE; 2019, {\tt arXiv:1812.08637}\par
\key Smith \paper Smith, R.; Nonlinear Kelvin and continental shelf waves; J. Fluid Mech.; 52 (1972) 379--391\par 
\key Smi2020a \paper Smith, D. A.; Revivals and Fractalization; Dyn. Sys. Web; 2020 (2020) 1--8\par
\key Sutherland \paper Sutherland, B.; Exact results for a quantum many-body problem in one-dimension. II; Phys. Rev. A; 5 (1972) 1372--1376\par
\key Talbot \paper Talbot, H.F.; Facts related to optical science. No. IV; Philos. Mag.; 9 (1836) 401--407\par
\key Thaller \book Thaller, B.; Visual Quantum Mechanics; Springer--Verlag, New York, 2000\par
\key VVS \paper Vrakking, M.J.J., Villeneuve, D.M., Stolow, A.; Observation of fractional revivals of a molecular wavepacket; Phys. Rev. A; 54  (1996) R37--40\par
\key Whitham \book Whitham, G.B.; Linear and Nonlinear Waves; John Wiley \& Sons, New York, 1974\par
\key YeaStr \paper Yeazell, J.A., Stroud, C.R.{, Jr.}; Observation of fractional revivals in the evolution of a Rydberg atomic wave packet; Phys. Rev. A; 43 (1991) 5153--5156\par
\key ZaZo \paper Zabrodin, A., Zotov, A.; Self-dual form of Ruijsenaars' Schneider models and ILW equation with discrete Laplacian; Nuclear Physics B; 927 (2018) 550--565\par

\endRefs

\page

\appendix
\section{Appendix}
\setcounter{equation}{0}\renewcommand{\theequation}{A.\arabic{equation}}

\subsection{Proof of \pr{propSC} 
} \label{sec:ProofOfPropSC}

\subsubsection{Definitions}

If $r>1$ then these series in equations~\eqref{eqn:SCdefn} converge by comparison with $1/n^2$.
When $r=1$, Dirichlet's test guarantees that these series converge for real $x$ not an integer multiple of $2\pi/k$.

We are particularly interested in the case $j\leq k$, under which condition both $S$ and $C$ can be expressed using a sum of polylogarithm functions, defined by the series
\Eq{polylog}
$$
    \Li_r(z) = \sum_{n=1}^\infty \frac{z^n}{n^r} = \sum_{n=0}^\infty \sum_{\gamma=1}^k \frac{z^{nk+\gamma}}{(nk+\gamma)^r}.
$$
For $x\in\mathbb{R}\setminus\set{2\pi l/k}{l\in\mathbb{Z}}$, we find
\begin{align*}
    \sum_{l=1}^k e^{-2\pi ilj/k} \Li_r\left(e^{2\pi \i l/k}e^{\i x}\right) &= \sum_{n=0}^\infty \sum_{\gamma=1}^k \frac{ e^{\i (nk+\gamma)x}  }{(nk+\gamma)^r}\ \sum_{l=1}^k e^{2\pi \i (\gamma-j)l/k}= k \sum_{n=0}^\infty \frac{e^{\i (nk+j)x}}{(nk+j)^r}.
\end{align*}
Thus,
\begin{subequations} \label{SC.polylogarithms}
\begin{align}
    S^k_{j,r}(x) &= \frac{1}{k} \Im\left[\; \sum_{l=1}^k e^{-2\pi \i jl/k} \Li_r\left(e^{2\pi i l/k}e^{ix}\right) \;\right], \\
    C^k_{j,r}(x) &= \frac{1}{k} \Re\left[ \;\sum_{l=1}^k e^{-2\pi \i jl/k} \Li_r\left(e^{2\pi i l/k}e^{ix}\right) \;\right].
\end{align}
\end{subequations}
If $r>1$ then this argument holds for all real $x$ hence equations~\eqref{SC.polylogarithms} are valid for all $r$ with $x\in\mathbb{R}\setminus\set{2\pi l/k}{l\in\mathbb{Z}}$ and also for all $r>1$ with $x\in\mathbb{R}$.

\subsubsection{Relations to other special functions}

It follows immediately from their definitions that, for $j,k,r\in\mathbb{N}^+$ with $j\leq k$,
\begin{align*}
    S^k_{j,r}(x) &= \Im\left[ \frac{e^{\i jx}}{k^r}L\left(\frac{kx}{2\pi},\frac{j}{k},r\right) \right] \\&= \Im\left.\left[ \frac{e^{\i jx}}{k^r}
    \prescript{}{1+r}{F}_r
    \left(
    \begin{matrix} 1,a,a,\ldots,a \\ a+1,a+1,\ldots,a+1 \end{matrix} ; e^{ikx}
    \right)
    \right]\right\rvert_{a=\frac{j}{k}}, \\
    C^k_{j,r}(x) &= \Re\left[ \frac{e^{\i jx}}{k^r}L\left(\frac{kx}{2\pi},\frac{j}{k},r\right) \right] \\&= \Re\left.\left[ \frac{e^{\i jx}}{k^r}
    \prescript{}{1+r}{F}_r
    \left(
    \begin{matrix} 1,a,a,\ldots,a \\ a+1,a+1,\ldots,a+1 \end{matrix} ; e^{ikx}
    \right)
    \right]\right\rvert_{a=\frac{j}{k}},
\end{align*}
where $L$ is the Lerch $\zeta$ function and $\prescript{}{1+r}{F}_r$ is a generalised hypergeometric function.

Let
\begin{equation}
    \Cl_r(x) = \begin{cases}\dsty \sum_{n=1}^\infty \frac{\cos(nx)}{n^{r}} ,\sstrut{20}\\\dsty \sum_{n=1}^\infty \frac{\sin(nx)}{n^{r}}, \end{cases} \qquad \Sl_r(x) = \begin{cases}\dsty \sum_{n=1}^\infty \frac{\sin(nx)}{n^{r}}, & \quad \mbox{if } r \mbox{ odd,}\sstrut{20} \\\dsty \sum_{n=1}^\infty \frac{\cos(nx)}{n^{r}}, & \quad\mbox{if } r \mbox{ even,} \end{cases}
\end{equation}
denote the Clausen functions $\Cl_r$ and the Glaisher-Clausen functions $\Sl_r$,~\cite{Lew1981a}.  In particular, for $r\in\{1,2\}$, these functions are given by
\begin{subequations} \label{eqn:SlCl12}
\begin{equation}
    \Sl_1(x) = \frac{\sign(x)\pi-x}{2}, \qquad \Sl_2(x) = \frac{x^2}{4} - \frac{\pi\lvert x \rvert}{2} + \frac{\pi^2}{6}, \qquad x\in(-\pi,\pi),
\end{equation}
and the Clausen functions are
\begin{equation}
    \Cl_1(x) = -\log\left\lvert2\sin\left(\frac{x}{2}\right)\right\rvert, \qquad \Cl_2(x) = \int_0^x \Cl_1(y)\dy, \qquad x\in(-\pi,\pi).
\end{equation}
\end{subequations}

Using the identities
\begin{equation}
    \Li_{r}(e^{\i x}) = \begin{cases} \Cl_{r}(x) + \i  \Sl_{r}(x) & \mbox{if } r \mbox{ odd,} \\ \Sl_{r}(x) + \i  \Cl_{r}(x) & \mbox{if } r \mbox{ even,} \end{cases}
\end{equation}
we obtain from equations~\eqref{SC.polylogarithms} the following representations:
\begin{subequations} \label{SC.Clausen}
\begin{align}
    S^k_{j,r}(x) &= \frac{1}{k}\sum_{l=1}^k \left[ -\sin\left(\frac{2\pi jl}{k}\right)\Cl_r\left(x+\frac{2\pi l}{k}\right) + \cos\left(\frac{2\pi jl}{k}\right)\Sl_r\left(x+\frac{2\pi l}{k}\right) \right], \\
    C^k_{j,r}(x) &= \frac{1}{k}\sum_{l=1}^k \left[ \cos\left(\frac{2\pi jl}{k}\right)\Cl_r\left(x+\frac{2\pi l}{k}\right) + \sin\left(\frac{2\pi jl}{k}\right)\Sl_r\left(x+\frac{2\pi l}{k}\right) \right],
\end{align}
for $r$ odd and, for $r$ even,
\begin{align}
    S^k_{j,r}(x) &= \frac{1}{k}\sum_{l=1}^k \left[ \cos\left(\frac{2\pi jl}{k}\right)\Cl_r\left(x+\frac{2\pi l}{k}\right) - \sin\left(\frac{2\pi jl}{k}\right)\Sl_r\left(x+\frac{2\pi l}{k}\right) \right], \\
    C^k_{j,r}(x) &= \frac{1}{k}\sum_{l=1}^k \left[ \sin\left(\frac{2\pi jl}{k}\right)\Cl_r\left(x+\frac{2\pi l}{k}\right) + \cos\left(\frac{2\pi jl}{k}\right)\Sl_r\left(x+\frac{2\pi l}{k}\right) \right].
\end{align}
\end{subequations}
Note that the Glaisher-Clausen functions are, on $(0,2\pi)$, polynomials closely related to the Bernoulli polynomials.
Therefore, on the interval $(-\pi,\pi)$, for $r$ even, $\Sl_r(x)$ is a polynomial in $\lvert x \rvert$ and, for $r$ odd, $\Sl_r(x)$ is the product of $\sign(x)$ with a polynomial in $\lvert x \rvert$.

\subsubsection{Elementary identities}

From the original definitions~\eqref{eqn:SCdefn}, for $M\in\mathbb{N}^+$,
\begin{equation} \label{eqn:TrigPolylog.Identities.Scale}
    M^r S^{Mk}_{Mj,r}\left(\frac{x}{M}\right) = S^k_{j,r}(x), \qquad 
    M^r C^{Mk}_{Mj,r}\left(\frac{x}{M}\right) = C^k_{j,r}(x),
\end{equation}
and
\begin{equation} \label{eqn:TrigPolylog.Identities.SumJ}
    \sum_{N=0}^{M-1} S^{Mk}_{j+Nk,r}(x) = S^k_{j,r}(x), \qquad 
    \sum_{N=0}^{M-1} C^{Mk}_{j+Nk,r}(x) = C^k_{j,r}(x).
\end{equation}
Equations~\eqref{eqn:TrigPolylog.Identities.Scale} follow from the original definitions~\eqref{eqn:SCdefn}.
Identities~\eqref{eqn:TrigPolylog.Identities.SumJ} follow from taking real and imaginary parts of
\begin{align*}
    \sum_{N=0}^{M-1} \frac{1}{Mk} \sum_{l=1}^{Mk} e^{-2\pi i (j+Nk) /Mk} \Li_r\left( e^{i(x+2\pi l/Mk)} \right)
    \hspace{-18em} & \\
    &= \frac{1}{Mk} \sum_{l=1}^{Mk} e^{-2\pi i lj /Mk} \sum_{N=0}^{M-1} e^{-2\pi i lN /M} \Li_r\left( e^{i(x+2\pi l/Mk)} \right) \\
    &= \frac{1}{Mk} \sum_{\substack{l=1:\\M\vert l}}^{Mk} e^{-2\pi i lj /Mk} M \Li_r\left( e^{i(x+2\pi l/Mk)} \right) 
    = \frac{1}{k} \sum_{l=1}^k e^{-2\pi i lj /k} \Li_r\left( e^{i(x+2\pi l/k)} \right).
\end{align*}

\subsubsection{Special cases}

If $k,j\in\{1,2\}$, then equations~\eqref{SC.Clausen} simplify to
\begin{align*}
    S^1_{1,r}(x) &= \begin{cases} \Sl_r(x) \\ \Cl_r(x) \end{cases} & S^2_{1,r}(x) &= \begin{cases} \frac{1}{2}\left(\Sl_r(x)-\Sl_r(x+\pi)\right) & \mbox{if } r \mbox{ odd,} \\ \frac{1}{2}\left(\Cl_r(x)-\Cl_r(x+\pi)\right) & \mbox{if } r \mbox{ even,} \end{cases} \\
    C^1_{1,r}(x) &= \begin{cases} \Cl_r(x) \\ \Sl_r(x) \end{cases} & C^2_{1,r}(x) &= \begin{cases} \frac{1}{2}\left(\Cl_r(x)-\Cl_r(x+\pi)\right) & \mbox{if } r \mbox{ odd,} \\ \frac{1}{2}\left(\Sl_r(x)-\Sl_r(x+\pi)\right) & \mbox{if } r \mbox{ even.} \end{cases}
\end{align*}
Identities~\eqref{eqn:TrigPolylog.Identities.Scale} may be used to deduce formulae for $S^2_{2,r}(x)$ and $C^2_{2,r}(x)$.

Much of this work concerns only the trigonometric polylogarithms in which $r=1$, which we will refer to as \emph{trigonometric hypergeometric functions}:
\begin{align*}
    S^k_{j}(x) &=  S^k_{j,1}(x) = \frac{1}{k} \Im\left[ \sum_{l=1}^k e^{-2\pi \i jl/k} \Li_1\left(e^{2\pi i l/k}e^{ix}\right) \right], \\
    C^k_j(x) &= C^k_{j,1}(x) = \frac{1}{k} \Re\left[ \sum_{l=1}^k e^{-2\pi \i jl/k} \Li_1\left(e^{2\pi i l/k}e^{ix}\right) \right],
\end{align*}
in which $\Li_1(z)=-\log(1-z)$.
In view of \eqref{eqn:SlCl12}, it follows that, on $(-\pi,\pi)$,
\Eq{SC112}
$$\eeq{    S^1_{1}(x) = \Sl_1(x) = \f2\bbk{\pi\sign(x)-x}, &C^1_{1}(x) = \Cl_1(x) = -\log\left\lvert2\sin\left(\frac{x}{2}\right)\right\rvert, \\
    S^2_{1}(x) = \f2\bbk{\Sl_1(x)-\Sl_1(x+\pi)}= \frac{\pi}{4}\,\sign(x), & C^2_{1}(x) = \f2\bbk{\Cl_1(x)-\Cl_1(x+\pi)}, \\
    S^1_{1,2}(x) = \Cl_2(x), & C^1_{1,2}(x) = \Sl_2(x).}
$$

\subsubsection{Smoothness, discontinuities and singularities} \label{ssec:SmoothDiscSing}

From their series representations in equations~\eqref{eqn:SCdefn}, it is immediate that 
$$
    S^k_{j,r},C^k_{j,r} \in H^{r-\frac{1}{2}-\epsilon}[-\pi,\pi] \mbox{ for } \epsilon>0.
$$
The smoothness described in \pr{propSC} follows by Sobolev's lemma~\cite{BrezisPDEBook}. 

From equations~\eqref{eqn:SlCl12} we see that, at $x=2\pi l$, $\Sl_1$ has jump discontinuities of height $1$, and $\Cl_1$ has logarithmic cusp singularities.
At the same points, $\Sl_2$ has corners and $\Cl_2$ has points of infinite gradient, but the second order functions $\Sl_2$ and $\Cl_2$, are continuous.
The characterisation of discontinuities and singularities described in \pr{propSC} follows by equations~\eqref{SC.Clausen}.

\subsubsection{Derivatives}

The derivative chains~\eqref{eqn:SCDerivatives} follow from equations~\eqref{SC.Clausen} and the derivatives for Clausen and Glaisher-Clausen functions~\cite[equations~(7.11--12) and~(7.15--16)]{Lew1981a}
\begin{equation} \label{eqn:SLCLDerivatives}
    \Cl'_{r+1} = (-1)^{r+1} \Cl_r, \qquad\qquad \Sl'_{r+1} = (-1)^r \Sl_r.
\end{equation}
The Clausen and Glaisher-Clausen functions of order 1 have distributional derivatives according with equations~\eqref{eqn:SLCLDerivatives} where
\begin{equation} \label{eqn:SlCl1Derivatives}
    \Cl_0(x) = \frac{1}{2}\cot\left(\frac{x}{2}\right), \qquad\qquad \Sl_0(x) = \pi\dide(x)-\frac{1}{2},
\end{equation}
are the order $0$ Clausen and Glaisher-Clausen functions and $\dide$ is the periodic extension of the $[-\pi,\pi]$ restriction of the Dirac delta function to $\mathbb{R}$.

In particular, for $j=1,2,\ldots,k-1$,
\begin{subequations} \label{eqn:dSC1}
\begin{align}
    \label{eqn:dS1.jgeneral}
    \frac{d}{dx}S^k_{j}(x) &= \frac{1}{2k}\sum_{l=1}^k\sin\left(\frac{2\pi jl}{k}\right)\cot\left(\frac{x}{2}+\frac{\pi l}{k}\right) + \frac{\pi}{k}\sum_{l=1}^k\cos\left(\frac{2\pi jl}{k}\right)\dide \left(x+\frac{2\pi l}{k}\right), \\
    \label{eqn:dS1.jk}
    \frac{d}{dx}S^k_{k}(x) &= \frac{\pi}{k}\sum_{l=1}^k\dide\left(x+\frac{2\pi l}{k}\right) - \frac{1}{2} = \Sl_0(kx), \\
\intertext{and, for $j=1,2,\ldots,k-1$}
    \label{eqn:dC1.jgeneral}
    \frac{d}{dx}C^k_{j,1}(x) &= -\frac{1}{2k}\sum_{l=1}^k\cos\left(\frac{2\pi jl}{k}\right)\cot\left(\frac{x}{2}+\frac{\pi l}{k}\right) + \frac{\pi}{k}\sum_{l=1}^k\sin\left(\frac{2\pi jl}{k}\right)\dide \left(x+\frac{2\pi l}{k}\right), \\
    \label{eqn:dC1.jk}
    \frac{d}{dx}C^k_{k,1}(x) &= \frac{-1}{2k}\sum_{l=1}^k\cot\left(\frac{x}{2}+\frac{\pi l}{k}\right) = -\Cl_0(kx).
\end{align}
\end{subequations}
Note that equation~\eqref{eqn:dC1.jk} can be seen as a special case of equation~\eqref{eqn:dC1.jgeneral}, but equation~\eqref{eqn:dS1.jk} is not a direct formulaic reduction of equation~\eqref{eqn:dS1.jgeneral}.

\subsection{Derivation of the expression \eq{smithif} for the convolution kernel of the Smith equation}

According to \Mathematica, the inverse Fourier transform of the Smith kernel
$$\shat_\delta(k) = \i\sqrt{k^2 + \fra \delta }\,$$
 is 
\Eq{sdx}
$$s_\delta(x) = - \i\sqrt{\frac 2{\pi\,\delta }}\;\frac{K_1\bpa{\abs x/\sqrt \delta \,}} {\abs x},$$
where $K_1(x)$ denotes the modified Bessel function of the second kind, \rf{DLMF; 10.25}

To derive \eq{sdx}, we start from  formula\fnote{There is a misprint in the formula as written.  The factor $\Gamma (\rho)$, which in our case is $\Gamma ( 3/2) = \sqrt \pi/2$, should be in the denominator of the right hand side.} 2.3.5.7 in \rf{PBM3; p. 323}:
\Eq{ift32}
$$\intoii \frac{e^{\i k\:x}\d k}{(k^2 + a^2)^{3/2}} 
= \frac2a\, \abs x \> K_1(a \,\abs x), \qquad a>0.$$
Next note that
$$\fhat(k) \equiv \odow k \sqrt{k^2 + a^2} = a^2 \,(k^2 + a^2)^{-3/2} .$$
According to \eq{ift32}, its  inverse Fourier transform equals
\Eq{FxK1} 
$$f(x) = \fra{\sqrt{2\pii}}\> \intoii{e^{\i k\:x} f(k) \d k} = \sqrt{\frac 2\pi} \> a\,\abs x \, K_1(a \,\abs x) .$$
However, 
$$\ghat(k) = \Into{-\:\infty}k\fhat(s)\d s = \frac k{\sqrt{k^2 + a^2}} + 1,$$
and so 
\Eq{hk}
$$\hhat(k) = \Into{-\:\infty}k\ghat(s) \d s = \sqrt{k^2 + a^2} + k.$$
On the other hand, we can express this double integration as a convolution
\Eq{hconv}
$$\hhat(k) = \intii s\rhat(k-s)\,\fhat(s) ,$$
where
$$\rhat(k) = \frac{k + \abs k}2 = \xcases{k, & k > 0,\\0,& k \leq 0,}$$
is the ramp function in the transform variable.  We note that the inverse Fourier transform of $k$ is  $\i \sqrt{2\pii} \> \dide^{\,\prime}(x)\ostrut{11}7$, a multiple of the derivative of delta, while the inverse Fourier transform of $\abs k$ is $-\, \sqrt{2/\pi}\> x^{-2}$.  Therefore, the inverse Fourier transform of the ramp function is the distribution
$$r(x) = -\fra{\sqrt{2\pii}\> x^2} + \i \sqrt{\frac \pi 2}\> \dide^{\,\prime}(x).$$
Consequently, the inverse Fourier transform of the convolution \eq{hconv} is obtained by multiplication:
\begin{align*}
h(x) = \sqrt{2\pii} \> r(x) \, f(x) &= -\, \frac{f(x)}{x^2} + \i \pi\, \dide^{\,\prime}(x)\, f(x) \\
&= -\, \frac{f(x)}{x^2} + \i \pi\, \bbk{f(0)\dide^{\,\prime}(x)- f'(0) \dide(x)}.
\end{align*}
The asymptotic properties of the Bessel function $K_1(x)$ at $x = 0$, \rf{DLMF}, imply
\Eq{F00}
$$\qeq{f(0) = \sqrt{\frac 2\pi} \,,\\f'(0) = 0.}$$
Therefore, using \eqas{hk}{F00}
$$\eeq{h(x) = -\, \sqrt{\frac 2\pi} \> \frac{a\,K_1(a \,\abs x)}{\abs x}   + \i \sqrt{2\pii}\> \dide'(x).}$$
The last term is the inverse Fourier transform of $k$, and hence, referring back to \eq{hk}, the first term is the inverse Fourier transform of $\sqrt{k^2 + a^2}$.  Replacing $a$ by $1/\sqrt\delta $ completes the derivation of \eq{sdx}.

\end{document}